\documentclass[a4paper,12pt]{amsart}

\usepackage{pgf,tikz,pgfplots}
\usepackage{tikz-cd}
\usepackage[obeyspaces]{url}
\usepackage{svg}
\pgfplotsset{compat=1.15}
\usetikzlibrary{matrix,arrows}
\usepackage{graphicx}

\usepackage[mathscr]{euscript} 
\usepackage{amsmath}
\usepackage{amsthm}
\usepackage{physics}
\usepackage{amssymb}
\usepackage{mathtools}

\usepackage{framed}
\usepackage{fancyhdr} 
\usepackage{vmargin}
\usepackage{afterpage}
\usepackage{hyperref}
\hypersetup{
    colorlinks=true,
}
\newcommand{\newabstract}[1]{%
  \par\bigskip
  \csname otherlanguage*\endcsname{#1}%
  \csname captions#1\endcsname
  \item[\hskip\labelsep\scshape\abstractname.]
\vspace*{0.25cm}}

\newcommand{\R}{\mathbb{R}} 

\newcommand{\vertiii}[1]{{\left\vert\kern-0.25ex\left\vert\kern-0.25ex\left\vert #1 
    \right\vert\kern-0.25ex\right\vert\kern-0.25ex\right\vert}}
    
\newcommand{\co}{\mathrm{co}}    
\newcommand{\clco}{\overline{\mathrm{co}}} 
\newcommand{\spn}{\mathrm{span} \,} 
\newcommand{\cod}{\mathrm{codim}} 
\newcommand{\dist}{\mathrm{dist}} 

\theoremstyle{plain}
\newtheorem{thm}{Theorem}[section] 
\newtheorem{df}[thm]{Definition} 
\newtheorem{ej}[thm]{Example} 
\newtheorem{lema}[thm]{Lemma} 
\newtheorem{prop}[thm]{Proposition} 
\newtheorem{cor}[thm]{Corollary} 
\newtheorem{rem}[thm]{Remark} 

\numberwithin{equation}{section} 

\title[Operators with the Kato property]{Operators with the Kato property on Banach spaces}

\author{Mar Jim\'enez-sevilla}
\address{ M. Jiménez-Sevilla. \newline Instituto de Matem\'atica Interdisciplinar and Departamento de An\'alisis Matem\'atico y Matem\'atica Aplicada, Facultad de Ciencias Matem\'aticas, Universidad Complutense de Madrid, Madrid 28040, Spain} \email{marjim@ucm.es}

\author{Sebasti\'an Lajara}
\address{S. Lajara. \newline Departamento de Matem\'aticas,  E.T.S.  de Ingenier\'ia Industrial, Universidad de Castilla-La Mancha,  02071 Albacete, Spain} \email{sebastian.lajara@uclm.es}

\author{Miguel Ángel Ruiz-Risue\~{n}o}
\address{M. \'A. Ruiz-Risue\~{n}o. \newline Departamento de Matem\'aticas y Ciencias de Datos, Facultad de Ciencias Econ\'omicas y Empresariales, Universidad CEU San Pablo, Madrid 28003, Spain} \email{miguelangel.ruizrisueno@ceu.es}
\subjclass[2010]{47A05,  46B20.}
\date{\today}
\keywords{Banach space,  operator with the Kato property,  separable quotient,  proper dense operator range, quasicomplemented subspace.}

\begin{document}

\baselineskip=15.5pt

\maketitle

\begin{abstract}
We consider a class of bounded linear operators between Banach spaces, which we call operators with the Kato property,  that includes the family of strictly singular operators between those spaces.  We show that if $T:E\to F$ is a dense-range operator with that property and $E$ has a separable quotient,  then for each proper dense operator range $R\subset E$ there exists a closed subspace $X\subset E$ such that $E/X$ is separable,   $T(X)$ is dense in $F$ and $R+X$ is infinite-codimensional.  If  $E^*$ is weak$^*$-separable,  the subspace $X$ can be built so that,  in addition to the former properties,  $R\cap X = \{0\}$.  Some applications to the geometry of Banach spaces are given.  In particular,  we provide the next extensions of well-known results of Johnson and Plichko: if $X$ and $Y$ are quasicomplemented but not complemented subspaces of a Banach space $E$ and $X$ has a separable quotient,  then $X$ contains a closed subspace $X_1\subset X$ such that $\dim (X/X_1)= \infty$ and $X_1$ is a quasicomplement of $Y$,  and that if $T:E\to F$ is an operator with non-closed range and $E$ has a separable quotient,  then there exists a weak$^*$-closed subspace $Z\subset E^*$ such that $T^*(F^*)\cap Z = \{0\}$.  Some refinements of these results, in the case that $E^*$ is weak$^*$-separable,  are also given.  Finally,  we show that if $E$ is a Banach space with a separable quotient,  then $E^*$ is weak$^*$-separable if,  and only if,  for every closed subspace $X\subset E$ and every proper dense operator range $R\subset E$ there exists a quasicomplement $Y$ of $X$ in $E$ such that $Y\cap R = \{0\}$.  
\end{abstract}

\section{Introduction}

A well-known result of Kato \cite{Kato} (c. f.  \cite[Proposition 2.c.4]{LT}) asserts that if $E$ and $F$ are Banach spaces and $T:E\to F$ is a  bounded linear operator whose restriction to every closed finite-codimensional subspace $Z\subset E$ is not an isomorphic embedding, then for every $\varepsilon> 0$ there exists an infinite-dimensional closed subspace $E_1\subset E$ such that $T|_{E_1}:E_1\to F$ (the restriction of $T$ to $E_1$) is a compact operator with $\|T|_{E_1}\|< \varepsilon$.  Observe that the above condition regarding finite-codimensional subspaces,  which will be referred to as the \textbf{Kato property} for the operator $T$, is satisfied if $T$ is strictly singular (that is,  the restriction of $T$ to every infinite-dimensional closed subspace $E_1\subset E$ is not an isomorphic embedding). 

A refinement of the former result,  in the case that $E$ is separable,  was obtained by  Fonf and Shevchik \cite{FS},  who proved that if $T:E\to F$ has the Kato property, then  the subspace $E_1$ can be built so that,  in addition to the features above, satisfies that
$\overline{T(E_1)} = \overline{T(E)}.$ More recently,  it has been shown in \cite{FLTZ} that,  under a strengthening of the Kato condition for the operator $T:E\to F$ and assuming also that $E$ is separable,  there exists a closed subspace $E_1\subset E$ which,  apart from the  property $\overline{T(E_1)} = \overline{T(E)}$,  has a special disposition with respect to a given closed subspace $L\subset E$ and an operator range $R$ in $E$ such that $R\cap L = \{0\}$ and ${\rm codim}\,  (R+L)=\infty$.  Recall that a vector subspace $R$ of a Banach space $E$ is said to be an \textbf{operator range} in that space if there exist a Banach space $G$ and an operator $A:G\to E$ such that $A(G)=R$.  
We notice that the class of operator ranges of a Banach space contains the one of its closed subspaces and is stable under finite sums and intersections (see e.g. \cite{Cross},  where we refer for the basic theory on operator ranges in Banach spaces).   The aforementioned result in \cite{FLTZ} reads as follows.  
\begin{thm}{\cite[Theorem 2.2]{FLTZ}}\label{main-FLTZ}
Let $E$ be a separable Banach space,  let $T:E\to F$ be an operator from $E$ into another Banach space $F$,  and let $L\subset E$ be a closed subspace of $E$ and $R$ be an operator range in $E$ such that $${\rm codim}\, (R+L) = \infty \quad \text{and} \quad R\cap L = \{0\}.$$ Assume that,  for every closed finite-codimensional subspace $Z\subset E$ containing $L$ there is a vector $x \in Z\setminus L$ such that $$\|T(x)\|\leq  \varepsilon \|Q_L(x)\|.$$ Then, there exists a closed subspace $E_1\subset E$ with the following properties: 
\begin{enumerate}
\item[(a)] $L\subset E_1$ and  $\dim \, (E_1/L)=\infty, $ 
\item[(b)] $E_1\cap R = \{0\}$,  and 
\item[(c)] $\overline{T(E_1)} = \overline{T(E)}.$
\end{enumerate}
\end{thm}

As an application of the former theorem,  it is proved in  \cite{FLTZ}  that if $H$ is a separable Hilbert space,  then for each closed subspace $H_0\subset H$ and each operator range $R$ in $H$ such that ${\rm codim}\,  (R+H_0)= \infty$ and $R\cap H_0 = \{0\}$ there exist closed subspaces $V,W\subset H$ such that $H= V\oplus_{\perp} W$, 
$H_0\subset V$ and $V\cap R = W\cap R = \{0\}$.  Some other results,  concerning the existence of quasicomplemented subspaces of a Banach space with a prescribed behaviour with respect to an operator range in that space,  have also been given in \cite{FLTZ} (recall that two closed subspaces of a Banach space $E$ are \textbf{quasicomplementary} if $X+Y$ is dense in $E$ and $X\cap Y = \{0\}$, in which case the subspace $Y$ is called a \textbf{quasicomplement} of $X$).   Corollaries 2.3 and 2.4 in \cite{FLTZ} provide the next extensions of well-known results of Murray-Mackey \cite{Mur, Mac} and James \cite[Theorem 3]{James} on quasicomplements: $(i)$  if $E$ is a separable Banach space then for every closed subspace $X\subset E$ and every infinite-codimensional operator range $R$ in $E$ containing $X$ there exists a closed subspace $Y\subset E$ such that  $Y$ is a quasicomplement of $X$ satisfying $Y\cap R = \{0\}$;  $(ii)$ if $X$ is a closed separable subspace of a Banach space $E$ and $Y\subset E$ is a proper quasicomplement of $X$ (that is, if $Y$ is not a complement of $X$),  then for every infinite-codimensional operator range $R$ in $X$ there exists a closed subspace $X_1\subset X$ such that $X_1\cap R = \{0\}$ and $Y$ is a quasicomplement of $X_1$. 

A relevant role in the Geometry of a Banach space $E$ is played by the class of proper dense operator ranges in that space (in the sequel, we denote that class with the symbol $\mathcal R_d(E)$).  Bennett and Kalton \cite{BK},  proved that a dense vector subspace $V$ of a Banach space $E$ is non-barrelled if, and only if, there exists $R \in \mathcal R_d (E)$ such that $V \subset R$.  On the other hand,  Saxon and Wilansky \cite{SW} showed that a Banach space $E$  contains a proper dense operator range if,  and only if,  $E$ has a separable (infinite-dimensional) quotient.   As remarked by Rosenthal \cite{R} (c. f.  \cite[Theorem 2.5]{Mujica}), the latter is equivalent to the existence of a closed separable quasicomplemented subspace of $E$. 
 Furthermore,  proper dense operator ranges possess nice structural properties.  For instance, Plichko \cite{P81} proved that if $E$ is a Banach space and  $R \in \mathcal R_d (E)$,  then the set $(E\setminus R) \cup \{0\}$ is \textbf{spaceable},  that is,  there exists a closed infinite-dimensional subspace $L \subset E$ such that $R \cap L = \{0\}$.   For more information on operator ranges and their applications we refer to the works \cite{COS} and \cite{FW}. 

In this paper, we provide some extensions of Theorem \ref{main-FLTZ} whenever $E$ is a Banach space with a separable quotient, and  obtain several applications to the study of quasicomplemented subspaces and the structure of operator ranges.  The work is divided into three sections,  apart from this Introduction.  In the next one,  we study a class of operators between Banach spaces which satisfy a Kato-like condition,  which is motivated by Therorem \ref{main-FLTZ} and  will be useful in the proofs of the results of the rest of the paper. 

In section 3,  we establish a partial extension of Theorem \ref{main-FLTZ} in the non-separable setting, proving that if $E$ is a Banach space with a separable quotient and $R\in \mathcal R_d(E)$, then under the hypothesis of that result there exists a closed subspace $E_1\subset E$ satisfying properties $(a)$ and $(c)$ and, instead of $(b)$, ${\rm codim} (R+E_1) = \infty$ (Theorem \ref{th:main}). In addition $E/E_1$ is separable and infinite-dimensional.  It is worth to mention that if the space is $E$ non-separable,  it is not possible in general to achieve the disjointness property $R \cap E_1 = \{0\}$ in Theorem \ref{main-FLTZ}  if we require that $\overline{T(E_1)} = \overline{T(E)}$ or $E/E_1$ is separable.  As a main application of our result, we deduce that if $X$ and $Y$ are proper quasicomplementary subspaces of a Banach space $E$ and $X$ has a separable quotient,  then there exists a closed subspace $X_1\subset X$ such that $\dim \, (X/X_1)=\infty$ and $X_1+Y$ is dense in $E$ (Corollary \ref{cr:james}).  This provides a generalization of a well-known theorem of Johnson \cite{Johnson}  which asserts that this statement holds true whenever $X$ is weakly compactly generated.  Another consequence of Theorem \ref{th:main}  yields that if $T : E \rightarrow F$ is a non surjective operator with dense range between Banach spaces and $E$ has a separable quotient, then there exists a weak$^*$-closed infinite-dimensional subspace $Z \subset E^*$ such that $T^* (F^*) \cap Z = \{0\}$,  being $T^*$ the adjoint operator of $T$ (Corollary \ref{cr:plichko}). This fact refines a result of Plichko \cite{P81},  who proved it in the case that the space $E$ is separable.

In the last section,  we provide a full generalization of Theorem \ref{main-FLTZ} in the case that the quotient $E/L$ has weak$^*$-separable dual and $R\in \mathcal R_d(E)$ (Theorem \ref{th:main-ws}).  This theorem leads to some extensions of the aforementioned statements from \cite{FLTZ} on quasicomplemented subspaces.  In particular,  we deduce that if $E$ is a Banach space with a separable quotient,  then the dual of $E$ is weak$^*$-separable if, and only if, for every closed subspace $X\subset E$ and every $R\in \mathcal R_d(E)$ containing $X$ there exists a closed subspace $Y\subset E$ such that $Y$ is a quasicomplement of $X$ and $R\cap Y = \{0\}$ (Corollary \ref{main2}).  Finally,  we obtain in Corollary \ref{plichko-dual-2} a refinement of Plichko's results on spaceability.

Our notation is standard, and can be found for instance in \cite{Hajek}.  We work with real normed spaces.  If $E$ is a normed space,  we denote by $B_E$ and $S_E$ the closed unit ball and the unit sphere  of that space,  respectively,  and by $I_E$ the identity operator on $E$. Unless otherwise stated, by ``\textbf{closed subspace}'', ``\textbf{operator}'' and ``\textbf{separable quotient}'' we mean ``\textbf{closed infinite-codimensional subspace}'', ``\textbf{bounded linear operator}'' and ``\textbf{infinite-dimensional separable quotient}'', respectively. Given a closed subspace $X$ of a Banach space $E$,  the symbol $Q_X$ stands for canonical quotient map $E \to E/X$.  If $A$ is a subset of $E$,  we denote by ${\rm span}\, (A)$ (respectively $[A]$) the linear subspace (respectively the closed linear subspace) generated by $A$,  by ${\rm co}\, (A)$ (respectively $\clco \, (A)$) the convex hull (respectively the closed convex hull) of $A$,  and  by $A^{\perp}$  its annihilator subspace,  that is,  $A^{\perp} = \{f\in E^*: \,  f(x)=0 \, \, \text{for all} \, x\in A\}$. 
Analogously,  if $B\subset E^*$,  we write $B_{\perp} = B^{\perp}\cap \pi(E)$,  being $\pi:E\to E^{**}$ the canonical isometry. If $V$ is a linear subspace of the dual of a Banach space $E$, then a sequence $\{u_n\}_n \subset E$ is called $V$-minimal if there exists a sequence $\{f_n\}_n \subset V$ biorthogonal to $\{u_n\}_n$. If $V = E^*$, we say that $\{u_n\}_n$ is minimal.

\section{Operators with the Kato property}
In this section we introduce the following extension of the Kato property,  motivated by Theorem \ref{main-FLTZ}, and obtain some results regarding this concept.
\begin{df}
	Let $T : E \rightarrow F$ be an operator between Banach spaces and $L$ be a closed subspace of $E$. We say that $T$ has the Kato property for $L$ if for every closed finite-codimensional subspace $W \subset E$ such that $L \subset W$ and every $\varepsilon > 0$ there exists $u \in W \setminus L$ such that
		$$\|T (u)\| \leq \varepsilon \|Q_L (u)\|.$$
\end{df}

As we mentioned in the Introduction, every strictly singular operator has the Kato property (in fact, if $T$ is strictly singular, then $T$ has the Kato property for every closed subspace $L \subset E$). However, the class of operators $T : E \rightarrow F$ with the Kato property is strictly larger than the one of strictly singular operators between $E$ and $F$.  Actually, as the following example shows, that class is not closed under finite sums.

\begin{ej}\normalfont
	Let $\{e_n\}_n$ be the canonical basis of $\ell_2$ and let $T : \ell_2 \rightarrow  \ell_2$ be the endomorphism given by $T (e_n) = 3^{-n} e_n$ for each $n \geq 1$. Now,  consider the endomorphisms $T_1, T_2 : \ell_2 \rightarrow \ell_2$ defined as 
$$T_1 (e_n) = \begin{cases} e_n & \text{if} \ n \ \text{is odd} \\
3^{-n} e_n & \text{if} \ n \ \text{is even}, \end{cases}
\quad \text{and} \quad T_2 (e_n) = \begin{cases} 3^{-n} e_n & \text{if} \ n \ \text{is odd} \\
e_n & \text{if} \ n \ \text{is even}. \end{cases}$$
It is clear that $T_1$ and $T_2$ have the Kato property and that $T_1 + T_2 = I_{\ell_2} + T$.  On the other hand,  $\|T\| \leq 1/2$, therefore $T_1 + T_2$ is an isomorphism on $\ell_2$.
\end{ej}

Next, we provide a characterization of the Kato property for subspaces.  Given an operator $T:E\to F$ between Banach spaces and a closed subspace $L\subset E$,  we write,  for each $u\in E$, 
$$\widetilde{T}_L \left(Q_L (u)\right) = Q_{\overline{T(L)}} \left( T (u)\right), \quad u \in E.$$
Notice that this expression is well-defined and yields a bounded linear operator $\widetilde{T}_L:E/L\to F/\overline{T(L)}.$

\begin{prop} \label{pr:kato}
	Let $T : E \rightarrow F$ be an operator between Banach spaces and $L$ be a closed subspace of $E$. Then the following assertions are equivalent:
	\begin{itemize}
		\item[(1)] $T$ has the Kato property for $L$.
		\item[(2)] The operator $\widetilde{T}_L$ has the Kato property. 
		\item[(3)] For every sequence $\{\varepsilon_n\}_n \subset (0, \infty)$ with $\varepsilon_n \rightarrow 0$ there is a sequence $\{u_n\}_n \subset E \setminus L$ such that $\{Q_L (u_n)\}_n$ is a normalized basic sequence in $E/L$ and
		$$\|T (u_n)\| \leq \varepsilon_n \quad \mathrm{\textit{for every}} \quad n \geq 1.$$
		\item[(3')] There exist a sequence $\{\varepsilon_n\}_n \subset (0, \infty)$ with $\varepsilon_n \rightarrow 0$ and a sequence $\{u_n\}_n \subset E \setminus L$ such that $\{Q_L (u_n)\}_n$ is a normalized basic sequence in $E/L$ and
		$$\|T (u_n)\| \leq \varepsilon_n \quad \mathrm{\textit{for every}} \quad n \geq 1.$$
	\end{itemize}
\end{prop}

\begin{proof}
	The implication $(3) \Rightarrow (3')$ is obvious.   We shall prove $(1) \Rightarrow (2)\Rightarrow (3)$ and $(3') \Rightarrow (1)$.
	
	$(1) \Rightarrow (2)$. Let $\widetilde{W} \subset E/L$ be a closed finite-codimensional subspace,  and set $W = Q_L^{-1} (\widetilde{W})$. Notice that $W$ has finite codimension in $E$, hence $W \setminus L \neq \emptyset$. Since $T$ has the Kato property for $L$,  there exists $u \in W \setminus L$ such that $\|T (u)\| \leq \varepsilon \|Q_L (u)\|$. Therefore,
	$$\|(\widetilde{T}_L \circ Q_L) (u)\| = \|(Q_{\overline{T(L)}} \circ T) (u)\| \leq \|T (u)\| \leq \varepsilon \|Q_L (u)\|.$$
	$(2) \Rightarrow (3)$. Let $\{\varepsilon_n\}_n \subset (0, \infty)$ be a sequence such that $\varepsilon_n \rightarrow 0$. Suppose that the operator $\widetilde{T}_L : E/L \rightarrow F/\overline{T(L)}$ has the Kato property. Then, arguing as in the proof of \cite[Proposition 2.c.4]{LT}, there exists a sequence $\{w_n\}_n \subset E$ such that $\{Q_L (w_n)\}_n$ is normalized basic in $E/L$ and such that $\|\widetilde{T}_L \circ Q_L (w_n)\| \leq 2^{-1} \varepsilon_n$ for each $n \geq 1$.  On the other hand, for each $w \in E$ we have 
	\begin{align*}
		\|(Q_{\overline{T(L)}} \circ T) (w)\| & = \inf \{ \|T(w) + z\| : z \in \overline{T(L)}\} = \inf \{ \|T(w) + z\| : z \in T(L)\} 
		\\ & = \inf \{ \|T (w + v)\| : v \in L\}.
	\end{align*}
	Thus,  there exists a sequence $\{v_n\}_n \subset L$ such that 
	$$\begin{aligned} \|T (v_n + w_n)\| & \leq \|(Q_{\overline{T(L)}} \circ T) (w_n)\| + 2^{-1} \varepsilon_n \|Q_L (w_n)\| \leq \varepsilon_n \|Q_L (w_n)\|\\ & = \varepsilon_n \|Q_L (v_n + w_n)\|      \end{aligned}$$
	for each $n \geq 1$.  Notice that if $u_n = v_n + w_n$ then $\{Q_L (u_n)\}_n$ is a normalized basic sequence in $E/L$.  Hence $\{u_n\}_n$ satisfies $(3)$. \par
$(3') \Rightarrow (1)$. Let $\{\varepsilon_n\}_n \subset (0, \infty)$ be a sequence such that $\varepsilon_n \rightarrow 0$ and $\{u_n\}_n \subset E \setminus L$ be a sequence such that $\{Q_L (u_n)\}_n$ is a normalized basic sequence in $E/L$ and $\|T (u_n)\| \leq \varepsilon_n$ for each $n \geq 1$. Let $\{f_n\}_n \subset L^{\perp}$ be the biorthogonal sequence of $\{u_n\}_n$,  and set $\alpha = \sup_n \|f_n\|$ (notice that $\alpha < \infty$ by the fact that $\{Q_L (u_n)\}_n$ is normalized basic).
	Let $W \subset E$ be a closed finite-codimensional subspace such that $L \subset W$,  and pick $\varepsilon> 0$.   We shall show the existence of a vector
$w \in (\spn \{u_n : n \geq 1\} \cap W) \setminus L$
such that $$\|T (w)\| \leq \varepsilon \|Q_L (w)\|.$$
Set $m = \cod_E \, (W) + 1$,  and take $k \geq 1$ such that $\varepsilon_n \leq \alpha^{-1} m^{-1} \varepsilon$ for each $n \geq k +1$.  Let us write
$$Z_{k,m} = [u_{k+1}, ..., u_{k+m}].$$
Since the sequence  $\{u_n\}_n$ is minimal we have $\dim \, Z_{k,m} = m$.  Bearing in mind that $W$ has codimension $m-1$ it follows that $\dim \, (W \cap Z_{k,m})\geq 1$.  Pick a non-zero vector $w \in W \cap Z_{k,m}$.  Then, there exist $\lambda_1, ..., \lambda_m \in \R$ such that 
$$w = \lambda_1 u_{k + 1} + \ldots + \lambda_m u_{k+m}.$$
As
$\|Q_L (w)\| = \sup \{f (w) : f \in S_{L^{\perp}}\}$
and $f_{k + i} (w) = \lambda_i$ for each $1 \leq i \leq m$ we get
$$\|Q_L (w)\| \geq \max_{1 \leq i \leq m} \|f_{k+i}\|^{-1} |f_{k+i} (w)|\geq \alpha^{-1} \max_{1 \leq i \leq m} |\lambda_i|.$$
On the other hand,
$$\|T (w)\| \leq \|T (\lambda_1 u_{k+1})\| + \ldots + \|T (\lambda_m u_{k+m})\| \leq m \max_{1 \leq i \leq m} \varepsilon_{k+i} \max_{1 \leq i \leq m} |\lambda_i|.$$
Therefore,
$$\|T (w)\| \leq \alpha m \|Q_L(w)\| \max_{1 \leq i \leq m} \varepsilon_{k+i} \leq \varepsilon \|Q_L (w)\|.$$
Consequently,  the operator $T$ has the Kato property for $L$.
\end{proof}

As application of the former proposition we obtain the following results, which provide some examples of operators with the Kato property.

\begin{cor}\label{cr:kato1}
	Let $E$ and $F$ be Banach spaces, $T : E \to F$ be an operator with non-closed range and $M \subset T (E)$ be a closed subspace. Then $T$ has the Kato property for $T^{-1} (M)$.
\end{cor}

\begin{proof}
	Set $L = T^{-1} (M)$. Suppose that $T$ does not have the Kato property for $L$.  Thanks to Proposition \ref{pr:kato} it follows that the induced operator $\widetilde{T}_L : E/L \rightarrow F/M$ does not enjoy the Kato property. Then there exists a closed finite-codimensional subspace $\widetilde{W} \subset E/L$ such that $\widetilde{T}_L|_{\widetilde{W}} : \widetilde{W} \rightarrow F/M$ is an isomorphic embedding. Since $\cod_{E/L} (\widetilde{W}) < \infty$ there exists a complement $\widetilde{N}$ of $\widetilde{W}$ in $E/L$,  hence
	$$\widetilde{T}_L (E/L) = \widetilde{T}_L (\widetilde{W}) + \widetilde{T}_L (\widetilde{N}),$$ and bearing in mind that $\widetilde{T}_L(\widetilde{N})$ is closed it follows that
$\widetilde{T}_L$ has closed range.  As $T (E) = Q_L^{-1} (\widetilde{T}_L (E/L))$ we have that $T(E)$ is closed,  which is a contradiction.
\end{proof}

\begin{cor}\label{cr:katoq1}
	Let $E$ be a Banach space and let $X, Y \subset E$ be closed subspaces. Then the following assertions are equivalent:
	\begin{itemize}
		\item[(1)] $X + Y$ is not closed in $E$.
		\item[(2)] The restriction to $X$ of the canonical quotient map $Q_Y : E \rightarrow E/Y$ has the Kato property for $X \cap Y$.
		\item[(3)] The restriction to $Y$ of the canonical quotient map $Q_X: E \rightarrow E/X$ has the Kato property for $X \cap Y$.
	\end{itemize}
\end{cor}

\begin{proof}  We only need to show the equivalence between $(1)$ and $(2)$.   Let $T$ be the restriction to $Y$ of the quotient map $Q_X:E\to E/X$.  If $X+Y$ is not closed then $T(Y)\neq E/X$,  and according to the previous corollary it follows that $T$ has the Kato property for the subspace $L = \ker T = X\cap Y$,  so $(1)\Rightarrow (2)$.  Conversely,  suppose that $T$ has the Kato property for $X\cap Y$,  then thanks to Proposition \ref{pr:kato} we deduce that 
the operator $\widetilde{T}:X/X\cap Y \to E/Y$ has the Kato property,  and because of its injectivity it follows that $\widetilde{T}(X/X\cap Y)$ is not closed in $E/Y$.  Therefore,  $X+Y$ is not closed in $E$.     
\end{proof}

\section{The main result}

In this section we provide an extension of Theorem \ref{main-FLTZ} assuming that the space $E$  has a  separable quotient.  We stress that the class of Banach spaces with a separable quotient is quite large: it includes $C(K)$ spaces,  being $K$ a Hausdorff compact space \cite{Lacey}, spaces with the separable complementation property (see \cite[p. 105]{Hajek}), and dual spaces \cite{Argyros}.  Nevertheless,  the classical problem asking for the existence of a separable quotient for every Banach space, still remains open.   The main result of this work reads as follows.

\begin{thm} \label{th:main}
	Let $E$ and $F$ be Banach spaces, $L \subset E$ be a closed subspace,  let $T : E \rightarrow F$ be an operator with the Kato property for $L$,  and $R \in \mathcal R_d (E)$.  If $\cod_E \, (L + R) = \infty$,  then there exists a closed subspace $E_1 \subset E$ such that 
\begin{enumerate}
\item[(a)] $L \subset E_1$, 
\item[(b)] $E/E_1$ is separable and infinite-dimensional, 
\item[(c)] $\cod_E \, (E_1 + R) = \infty$, and
\item[(d)] $\overline{T (E_1)} = \overline{T (E)}.$
\end{enumerate}
\end{thm}

As we mentioned in the Introduction,  Bennett and Kalton \cite{BK} proved that if $V$ is a dense vector subspace of a Banach space $E$,  then $V$ is non-barrelled if,  and only if,  there exists $R\in \mathcal R_d(E)$ such that $V\subset R$.  Thus,   the former theorem holds true if $R$ is a dense non-barrelled subspace of $E$.  We also point out that if $E$ is non-separable then,  under the hypothesis of this theorem (and assuming that $R\cap L = \{0\}$ instead that ${\rm codim}_E (R+L)=\infty$),  it is not possible in general to find a closed subspace $E_1\subset E$ satisfying  
$R\cap E_1 = \{0\}$ (instead of property $(c)$),  and such that either $\overline{T(E_1)} = \overline{T(E)}$ or 
$E/E_1$ is separable.  Indeed,  as remarked in \cite[p. 206]{FLTZ},  the classical theorem of Lindenstrauss \cite{L} (see also \cite[Corollary 5.89]{Hajek}) ensuring that $c_0(\Gamma)$ is not quasicomplemented in $\ell^{\infty}(\Gamma)$ whenever the set $\Gamma$ is uncountable yields that if $T:\ell^{\infty}(\Gamma)\to \ell^{\infty}(\Gamma)/c_0(\Gamma)$ is the canonical quotient map and $R$ is a proper dense operator range in $\ell^{\infty}(\Gamma)$ containing $c_0(\Gamma)$ then, for any closed subspace $E_1\subset \ell^{\infty}(\Gamma)$ such that $R\cap E_1 = \{0\}$ we have that $T(E_1)$ is not dense in $\ell^{\infty}(\Gamma)/c_0(\Gamma)=T(\ell^{\infty}(\Gamma))$ (notice that the operator $T$ satisfies the Kato property).  On the other hand,  in \cite[Theorem 3.2]{JL3} it is shown that if $E$ is a Banach space with the property that for every $R\in \mathcal R_d(E)$ there exists a closed subspace $E_1\subset E$ such that $E/E_1$ is separable and $R\cap E_1 = \{0\}$,  then $E^*$ is weak$^*$-separable.  In the next section,  we shall show that if the space $(E/L)^*$ is weak$^*$-separable and $L\cap R = \{0\}$ then,  under the hypothesis of Theorem \ref{th:main},  there exists a closed subspace $E_1\subset E$ satisfying properties $(a), (b),(d)$ and $E_1\cap R = \{0\}$. 

A key ingredient in the proof of Theorem \ref{th:main} is the following technical lemma.

\begin{lema}\label{lm:sw}
	Let $E$ and $F$ be Banach spaces, $T : E \rightarrow F$ be an operator, $M \subset E$ be a closed subspace,  and let $\{z_n, f_n\}_n \subset E \times M^{\perp}$ be a biorthogonal system.   Assume there exists a number $\alpha > 1$ satisfying:
	\begin{itemize}
		\item[(i)] $\|Q_M (z_n)\| \leq \alpha^n \|f_n\|^{-1}$ and
		\item[(ii)] $\|T (z_n)\| \leq \alpha^{-2n} \|Q_M (z_n)\|$ 
	\end{itemize}
	for each $n\geq 1$.  Then, for every sequence $\{\varepsilon_n\}_n \subset (0, \infty)$ and every sequence $\{W_n\}_n$ of closed finite-codimensional subspaces of $E$ there exists a block sequence $\{u_n\}_n$ of $\{z_n\}_n$ with the following properties: 
	\begin{itemize}
		\item[(1)] $\{Q_M (u_n)\}_n$ is a normalized basic sequence in $E/M$,
		\item[(2)] $u_n \in \bigcap_{i=1}^n W_i$ for each $n \geq 1$, and
		\item[(3)] $\|T (u_n)\| \leq \varepsilon_n$ for each $n \geq 1$.
	\end{itemize}
\end{lema}

\begin{proof}
	Consider, for every $k, n \geq 1$, the subspace
	$$Z_{k,n} = [z_{k+1}, ..., z_{k+n}].$$
	Notice that, since $\{z_n\}_n$ is $M^{\perp}$-minimal, $\dim \, Z_{k,n} = n$. We claim that for each $z \in Z_{k,n}$,
		\begin{equation*}
			\|T (z)\| \leq \alpha^{-k+n} n \|Q_M (z)\|.
		\end{equation*}
	First, we will show that $Q_M|_{Z_{k,n}}$ is injective. Pick $u \in Z_{k,n}$. Then there exist $\lambda_1, ..., \lambda_n \in \R$ such that
	 	$$u = \lambda_1 z_{k+1} + \ldots + \lambda_n z_{k+n}.$$
	Suppose that $Q_M (u) = 0$. As $\{f_n\}_n \subset M^{\perp}$ we get $f_j (u) = 0$ for each $j \geq 1$.  Hence, for each $1 \leq i \leq n$ we have 
	 	$$0 = f_{k+i} (u) = f_{k+i} (\lambda_1 z_{k+1} + \ldots + \lambda_n z_{k+n}) = \lambda_i f_{k+i} (z_{k+i}) = \lambda_i,$$
	Therefore, $u = 0$.\par
	Set
		$$\eta = \inf \left \{ \left \| \sum_{i=1}^n \lambda_i \frac{Q_M (z_{k+i})}{\|Q_M(z_{k+i})\|}  \right \| : \,  |\lambda_1| + \ldots + |\lambda_n| = 1 \right \}.$$
	We assert that $$\eta \geq \alpha^{-k-n} n^{-1}.$$ Indeed, for each $\lambda_1, ..., \lambda_n \in \R$ and each $i=1,\ldots,n$ we have
		\begin{align*}
		 	\left \| \sum_{j=1}^n \lambda_j  \frac{Q_M (z_{k+j})}{\|Q_M(z_{k+j})\|} \right \| & \geq \|f_{k+i}\|^{-1} \left | f_{k+i} \left( \sum_{j=1}^n \lambda_j \frac{Q_M (z_{k+j})}{\|Q_M(z_{k+j})\|}  \right) \right |  \\
		 	& = |\lambda_i| \|Q_M (z_{k+i})\|^{-1} \|f_{k+i}\|^{-1} \geq \alpha^{-k-i} |\lambda_i| \\
		 	& \geq \alpha^{-k-n} |\lambda_i|.
		\end{align*}
	Thus $\eta \geq \alpha^{-k-n} n^{-1}$.  Now,  pick $z \in Z_{k,n}$ such that $\|Q_M (z)\| = 1$. Then 
		$$Q_M (z) \in B_{E/M} \cap Q_M (Z_{k,n}) = B_{Q_M (Z_{k,n})}.$$
	Notice that
		$$\eta B_{Q_M(Z_{k,n})} \subset \co \{\pm \|Q_M (z_{k+i})\|^{-1} Q_M (z_{k+i}) : 1 \leq i \leq n\}.$$
	Since $Q_M (z) \in B_{Q_M (Z_{k,n})}$ and $Q_M|_{Z_{k,n}}$ is injective, we get
		$$z \in \eta^{-1} \co \{\pm \|Q_M (z_{k+i})\|^{-1} z_{k+i} : 1 \leq i \leq n\}.$$
	Therefore, by (ii),
		\begin{align*}
		 	\|T (z)\| & \leq \eta^{-1} \sup \big{\{} \|T (y)\| :\, y \in \co \{ \pm \|Q_M (z_{k+i})\|^{-1} z_{k+i} : 1 \leq i \leq n\} \big{\}}  \\
		 	& = \eta^{-1} \sup \{ \|Q_M (z_{k+i})\|^{-1} \|T (z_{k+i})\| :\, 1 \leq i \leq n\} \leq \alpha^{-2k} \eta^{-1} \\
		 	& \leq \alpha^{-k+n} n.
		\end{align*}
	Now,  we proceed with the construction of the sequence $\{u_n\}_n$.  We elaborate on the argument of Mazur's theorem (see \cite[Theorem 4.19]{Fabian}). Let $\{\delta_n\}_n \subset (0,1)$ be a sequence such that $\Pi_n (1 + \delta_n) < \infty$.   
	Set $\mu_1 = \cod_E \, (W_1)$ and  $q_1 = \mu_1 + 1$. Pick $p_1 \geq 1$ such that $\alpha^{-p_1+q_1} \leq q_1^{-1} \varepsilon_1$,  and consider the subspace
		$$Z_{p_1,q_1} = [z_{p_1+1}, ..., z_{p_1+q_1}].$$
	Notice that $\dim \, Z_{p_1,q_1} = q_1$. Therefore, there exists a vector
		$$u_1 \in W_1 \cap Z_{p_1,q_1}$$
	such that $\|Q_M (u_1)\| = 1$. Since $u_1 \in Z_{p_1,q_1}$, we obtain
	$$\|T (z_1)\| \leq \alpha^{-p_1+q_1} q_1 \leq \varepsilon_1.$$
	Fix $n \geq 1$,  and suppose that there exists a ``finite block sequence'' $u_1, ..., u_n \in \spn \{z_k : k \geq 1\}$ such that $\|Q_M (u_i)\| = 1$, $\|T( u_i)\| \leq \varepsilon_i$, and 
		$$u_i \in \bigcap_{j=1}^i W_j$$
	for each $1 \leq i \leq n$ and 
		$$\|y\| \leq (1 + \delta_i) \|y + \lambda Q_M (u_{i+1})\|$$
	whenever $1 \leq i \leq n-1$,  $y \in [Q_M (u_1), ..., Q_M (u_i)]$ and $\lambda \in \R$.  
	
	Since $u_1, ..., u_n \in \spn \{z_k : k \geq 1\}$ there exists $k_n \geq 1$ such that $u_1, ..., u_n \in [z_1, ..., z_{k_n}]$.  Let us write 
		$$U_n = [Q_M (u_1), ..., Q_M (u_n)],$$ and 
	let $\{\omega_{j,n}\}_{j=1}^{m_n}$ be a subset of $S_{M^{\perp}}$ such that $\{\omega_{j,n}|_{U_n}\}_{j=1}^{m_n}$ is a $(\delta_n/2)$-net for $S_{U_n^*}$.  Set $\mu_{n+1} = \sum_{i=1}^{n+1} \cod_E \, (W_i)$ and  $q_{n+1} = m_n + \mu_{n+1} + 1$.   Pick $p_{n+1} \geq k_n$ such that $\alpha^{-p_{n+1}+q_{n+1}} \leq q_{n+1}^{-1} \varepsilon_{n+1}$,  and consider the subspace
		$$Z_{p_{n+1},q_{n+1}} = [z_{p_{n+1}+1}, ..., z_{p_{n+1}+q_{n+1}}].$$
	Notice that $\dim \, Z_{p_{n+1},q_{n+1}} = q_{n+1}$. Therefore, since
	$$\cod_E \, \left(\left(\bigcap_{i=1}^{n+1} W_i\right) \cap \left(\bigcap_{j=1}^{m_n} \ker \omega_{j,n}\right)\right) \leq m_n + \mu_{n+1} = q_{n+1} - 1,$$
	there exists a vector
		$$u_{n+1} \in \left(\bigcap_{i=1}^{n+1} W_i\right) \cap \left(\bigcap_{j=1}^{m_n} \ker \omega_{j,n}\right) \cap Z_{p_{n+1},q_{n+1}}$$
	such that $\|Q_M (u_{n+1})\| = 1$.  Bearing in mind that $u_{n+1} \in \bigcap_{j=1}^{m_n} \ker \omega_{j,n}$ we get
		$$\|u\| \leq (1 + \delta_n) \|u + \lambda Q_M (u_{n+1})\|$$
	whenever $u \in [Q_M (u_1), ..., Q_M (u_n)]$ and $\lambda \in \R$.  As $u_{n+1} \in Z_{p_{n+1},q_{n+1}}$ we obtain
		$$\|T (u_{n+1})\| \leq \alpha^{-p_{n+1}+q_{n+1}} q_{n+1} \leq \varepsilon_{n+1}.$$
	We have proved that $\{Q_M (u_n)\}_n$ is a normalized basic sequence in $E/M$ such that $u_n \in \bigcap_{i=1}^n W_i$ and $\|T (u_n)\| \leq \varepsilon_n$ for each $n \geq 1$.
\end{proof}

As a consequence of the previous lemma, we obtain the following corollary.

\begin{cor}\label{cr:sw}
	Let $E$ and $F$ be Banach spaces,  let $T : E \rightarrow F$ be an operator and $M$ be a closed subspace of $E$.   If there exist a biorthogonal system $\{z_n, f_n\}_n \subset E \times M^{\perp}$ and a number $\alpha > 1$ such that 
	\begin{itemize}
		\item[(i)] $\|Q_M (z_n)\| \leq \alpha^n \|f_n\|^{-1}$ and 
		\item[(ii)] $\|T (z_n)\| \leq \alpha^{-2n} \|Q_M (z_n)\|$ 
	\end{itemize}
	for each $n\geq 1$,  then $T$ has the Kato property for $M$.
\end{cor}

\begin{proof}
	Let $\{\varepsilon_n\}_n \subset (0,\infty)$ be a sequence such that $\varepsilon_n \rightarrow 0$. The former lemma guarantees the existence of a sequence $\{u_n\}_n$ such that $\{Q_M (u_n)\}_n$ is a normalized basic sequence in $E/M$ and $\|T (u_n)\| \leq \varepsilon_n$ for each $n \geq 1$,  and thanks to Proposition \ref{pr:kato} we deduce that $T$ has Kato property for $M$.
\end{proof}

We shall also need the following result.

\begin{lema}\label{lm:min}
Let $V$ be a bounded closed convex symmetric subset of a Banach space $E$. Then the Minkowski's functional $\mu_V$ provides a norm $\|\cdot\|_V = \mu_V (\cdot)$ in $\spn V$ such that $(\spn V, \|\cdot\|_V)$ is a Banach space. In particular, if $V$ has empty interior and $\spn V$ is dense in $E$, then $\spn V \in \mathcal R_d (E)$.
\end{lema}

\begin{proof}
The first part of the lemma appears in \cite[pg. 73]{Fabian}. Suppose that $V$ has empty interior. Set $R = \spn V$. Notice that
$$R = \bigcup_m m V.$$
Then, by the Baire's category theorem, the subspace $R$ has also empty interior. Let $\psi : (R, \|\cdot\|_V) \rightarrow (E, \|\cdot\|)$ be the inclusion operator of $R$. Notice that since $V$ is bounded, $\psi$ is continuous. By the open mapping theorem,  as $R$ has empty interior, it is a proper vector subspace. On the other hand, since $R$ is dense, it cannot be a closed subspace. Thus,  $\cod_E \, (R) = \infty$,  and hence $R \in \mathcal R_d (E)$.
\end{proof}

\noindent \emph{Proof of Theorem \ref{th:main}. } We divide the proof into several steps.

	\textsc{Step 1}. There exist vector sequences $\{v_n\}_n, \{w_n\}_n, \{z_n\}_n \subset S_E$, a functional sequence $\{f_n\}_n \subset E^*$, and a sequence $\{B_n\}_{n \geq 0}$ of bounded closed convex symmetric subsets of $E$ such that, for every $n \geq 1$:
		\begin{itemize}
			\item[(1)] $f_n (w_n) = 1$ and $f_n (w_k) = 0$ if $k \neq n$,  $k \geq 1$,
			\item[(2)] $f_n (z_k) \neq 0$ if and only if $k = n$,
			\item[(3)] $\|T (z_n)\| \leq 3^{-2n} \|f_n\|^{-1} f_n (z_n)$,
			\item[(4)] $f_n (v_n) \geq 2^{-1} \|f_n\|$ and $f_k (v_n) = 0$ if $k > n$,
			\item[(5)] $B_n = B_{n-1} + \{ \xi v_n + \xi' w_n + \xi'' z_n : |\xi| + |\xi'| + |\xi''| \leq 1\}$, and
			\item[(6)] $\sup_{B_{n-1}} f_n \leq 2^{-n}$.
		\end{itemize}
		
		Let $G$ be a Banach space and $A : G \rightarrow E$ be an operator with $\|A\|=1$ such that $R = A (G)$, and set $B_0 = \overline{A (B_G)}$. Notice that $B_0$ is bounded, closed, convex and symmetric, has empty interior and $\spn B_0$ is dense in $E$.  Thus,  Lemma \ref{lm:min} yields that $\spn B_0 \in \mathcal R_d (E)$.\par
		Since $\cod_E \, (\spn B_0) = \infty$, there exists a vector $w_1 \in S_E \setminus \spn B_0$.   Thanks to the Hahn-Banach separation theorem, we can find a functional $g_1 \in E^*$ such that $g_1 (w_1) = 1$ and $\sup_{B_0} g_1 \leq 2^{-2}$.  Let $Y_1 \subset E$ be a complement of $[w_1]$. Since $T$ has the Kato property for $L$ we have that $T$ posesses the Kato property. This, combined with the fact that $\cod_E \, (Y_1) < \infty$,  guarantees the existence of a vector $z_1 \in S_{Y_1}$ such that
		\begin{equation}\label{th:fltz:eq:1}
			\|T (z_1)\| \leq 2^{-2} 3^{-2} (\|g_1\| + 1)^{-1} \dist ([w_1], S_{Y_1}).
		\end{equation}
		We may assume that $g_1 (z_1) \geq 0$. By the Hahn-Banach theorem, we can find $h_1 \in S_{E^*}$ such that $h_1 (w_1) = 0$ and $h_1 (z_1) = \dist (z_1, [w_1]) > 0$. Now, define
		$$f_1 = g_1 + 2^{-2} h_1.$$
		It is immediate that
		$f_1 (w_1) = 1,$
		$f_1 (z_1) \geq 2^{-2} \nu_1, $
		and $\sup_{B_0} f_1 \leq 2^{-1}.$
		Bearing in mind that $\|f_1\| \leq \|g_1\| + 2^{-2} \leq \|g_1\| + 1$ and $f_1 (z_1) \geq 2^{-2} \dist ([w_1], S_{Y_1})$,  thanks to \eqref{th:fltz:eq:1} we get 
		$$\|T (z_1)\| \leq 3^{-2} \|f_1\|^{-1} f_1 (z_1).$$
		Finally,  choose a vector $v_1 \in S_E$ such that $f_1 (v_1) \geq 2^{-1} \|f_1\|$. \par
		Now, pick $n \geq 1$ and suppose that there exist vectors $v_1, ..., v_n, w_1, ..., w_n, z_1, ..., z_n \in S_E$, functionals $f_1, ..., f_n \in E^*$ and sets $B_0, ..., B_{n-1} \subset E$ such that for every $1 \leq i, j  \leq n$:
		\begin{itemize}
			\item[(i)] $f_i (w_j) = \delta_{i,j}$,   
			\item[(ii)] $f_j (z_i) \neq 0$ if and only if $ j\neq i $,
			\item[(iii)] $\|T (z_i)\| \leq 3^{-2i} \|f_i\|^{-1} f_i (z_i)$,
			\item[(iv)] $f_i (v_i) \geq 2^{-1} \|f_i\|$, and $f_j (v_i) = 0$ whenever $j\neq i$,  
			\item[(v)] $B_i = B_{i-1} + \{ \xi v_i + \xi' w_i + \xi'' z_i : |\xi| + |\xi'| + |\xi''| \leq 1\}$ if $1 \leq i \leq n-1$, and
			\item[(vi)] $\sup_{B_{i-1}} f_i \leq 2^{-i}$.
		\end{itemize}
		Set
		$$B_n = B_{n-1} + \{ \xi v_n + \xi' w_n + \xi'' z_n : |\xi| + |\xi'| + |\xi''| \leq 1\}.$$
		Taking into account that $B_0$ is bounded, closed, convex, symmetric, has empty interior, and $\spn B_0$ is dense in $E$, the same occurs with $B_n$. Therefore, by Lemma \ref{lm:min}, $\spn B_n \in \mathcal R_d (E)$. In particular, $\cod_E \, (\spn B_n) = \infty$. Let us write
		$$N_n = [v_1, ..., v_n, w_1, ..., w_n, z_1, ..., z_n].$$
		Since $N_n \subset \spn B_n$ we have 
		$$Q_{N_n} (\spn B_n) \in \mathcal R_d (E/N_n).$$
		Set $\widetilde{B}_n = \overline{Q_{N_n} (B_n)}^{E/N_n}$.  As $\spn (Q_{N_n} (B_n)) = Q_{N_n} (\spn B_n) \in \mathcal R_d (E/N_n)$ we get
		$$\spn \widetilde{B}_n \in \mathcal R_d (E/N_n).$$ Therefore, 
		according to \cite[Proposition 1.10]{Dixmier} we obtain $Q_{N_n}^{-1} (\spn \widetilde{B}_n) \in \mathcal R_d (E)$. Let $G_n$ be a Banach space and $A_n : G_n \rightarrow E$ be an operator such that $\|A_n\| = 1$ and $A_n (G_n) = Q_{N_n}^{-1} (\spn \widetilde{B}_n)$,  and set $V_n = \overline{A_n (B_{G_n})}$. Notice that $V_n$ is bounded, closed, convex, symmetric, has empty interior and $\spn V_n$ is dense in $E$. Hence,  by Lemma \ref{lm:min}, $\spn V_n \in \mathcal R_d (E)$. In particular, $\cod_E \, (\spn V_n) = \infty$. Therefore, there exists a vector
			\begin{equation}\label{th:fltz:eq:2}
		  	w_{n+1} \in \left( S_E \cap \left( \bigcap_{i=1}^n \ker f_i \right) \right) \setminus \spn V_n.
			\end{equation}
		Observe that,
		since $\widetilde{B}_n \subset \spn (Q_{N_n} (V_n))$ and $N_n \subset \spn V_n$, we have  
		$$Q_{N_n} (w_{n+1}) \in (E/N_n) \setminus \spn \widetilde{B}_n.$$
		As $N_n^{\perp}$ identifies with $(E/N_n)^*$,  the Hahn-Banach separation theorem yields a functional $g_{n+1} \in N_n^{\perp}$ such that $g_{n+1} (w_{n+1}) = 1$ and $\sup_{\widetilde{B}_n} g_{n+1} \leq 2^{-n-2}$.\par
		Let $Y_{n+1} \subset E$ be a complement of $N_n + [w_{n+1}]$ in $E$. Since $T$ has the Kato property and $\cod_E \, (Y_{n+1}) < \infty$ there exists a norm-one vector
		\begin{equation}\label{th:fltz:eq:4}
		 	z_{n+1} \in Y_{n+1} \cap \left( \bigcap_{i=1}^n \ker f_i \right)
		\end{equation}
		satisfying
		\begin{equation}\label{th:fltz:eq:5}
			\|T (z_{n+1})\| \leq 2^{-n-2} 3^{-2(n+1)} (\|g_{n+1}\| + 1)^{-1} \dist (N_n + [w_{n+1}], S_{Y_{n+1}}).
		\end{equation}
		We may assume that $g_{n+1} (z_{n+1}) \geq 0$. Notice that $Q_{N_n} (z_{n+1}) \neq 0$. By the Hahn-Banach theorem, there exists $h_{n+1} \in S_{N_n^{\perp}}$ such that $h_{n+1} (w_{n+1}) = 0$ and $h_{n+1} (z_{n+1}) = \dist (z_{n+1}, N_n + [w_{n+1}]) > 0$. Now, define
		$$f_{n+1} = g_{n+1} + 2^{-n-2} h_{n+1},$$
		and pick $v_{n+1} \in S_E$ such that
		$$f_{n+1} (v_{n+1}) \geq 2^{-1} \|f_{n+1}\|.$$
		The construction is finished and it is obvious that the sequences $\{v_n\}_n$, $\{w_n\}_n$, $\{z_n\}_n$, $\{f_n\}_n$, and $\{B_n\}_{n \geq 0}$ satisfy properties (1) - (6).\par	
		
	\textsc{Step 2}: There exists a closed subspace $M \subset E$ such that $E/M$ is infinite-dimensional and separable,  $\cod (M + R) = \infty$ and $T$ has the Kato property for $M$.\par
		To prove this statement, consider sequences $\{v_n\}_n$, $\{w_n\}_n$, $\{z_n\}_n$, $\{f_n\}_n$, and $\{B_n\}_{n \geq 0}$ satisfying the properties from Step 1. Let us write
		$$M = \bigcap_n \ker f_n.$$ 
		
		The constructions in \cite[Proposition 1.10]{SW} and \cite[Proposition 3.3]{JL3} yield respectively that $E/M$ is infinite-dimensional and separable, and $\cod_E \, (M + R) = \infty$. We include the details for the sake of completeness. The fact that $M$ is infinite-codimensional in $E$ comes from property (1) from the previous step. Let us prove that $E/M$ is separable.  Let us write,  for each $n\geq 1$,  $M_n = \bigcap_{j \geq n} \ker f_j$ (notice that $M_1 = M$). We claim that $\bigcup_n M_n$ is dense in $E$.  Consider the map $\psi_n : B_{n-1} \rightarrow E$ defined as 
		$$\psi_n (u) = u - \sum_{j \geq n} f_j (u) w_j, \quad u \in B_{n-1}.$$
		Pick $n \geq 1$ and $u \in B_{n-1}$. Then, by property (1) from the previous step we have that $|f_j (u)| \leq 2^{-j}$ for each $j \geq n$, hence $\|\psi_n (u) - u\| \leq 2^{-n+1}$ for each $n \geq 1$. On the other hand,  using again property (1), for every $i \geq n$ we get
		$$f_i (\psi_n (u)) = f_i \left( u - \sum_{j \geq n} f_j (u) w_j \right) = f_i (u) - \sum_{j \geq n} f_j (u) f_i (w_j) = 0,$$
		hence $\psi_n (u) \in M_n$.  Since $\spn B_n$ is dense in $E$ for every $n \geq 0$ it follows that $\bigcup_n M_n$ is dense in $E$.\par
		As $\dim (M_{n+1}/M_n) = 1$ and $w_n \in M_{n+1} \setminus M_n$ for each $n \geq 1$ we infer that
		$$M_{n+1} = M \oplus [w_1, ..., w_n] \quad \text{for \ every} \quad n \geq 1.$$
		Let $g \in M^{\perp}$ be a functional such that $g (w_n) = 0$ for each $n \geq 1$. Then $g$ vanishes on $M_n$ for each $n \geq 1$. Since $\bigcup_n M_n$ is dense in $E$, it follows that $g = 0$. Thus, $\{Q_M (w_n)\}_n$ is linearly dense in $E/M$,  a consequently  $E/M$ is separable.\par
		Next, we shall show that $\cod (R + M) = \infty$. Consider the vector subspace
		$$S = \left\{u \in E \, : \, \text{the sequence} \quad \{n f_n (u)\}_n \quad \text{converges} \right\},$$
		and set,  for each $u\in S$, 
		$$\vertiii{u} = \|u\| +  \sup \left\{n |f_n (u)|: \, n\geq 1 \right\}.$$
		It is easy to check that $\vertiii{\cdot}$ is a complete norm on $S$.  Since $\vertiii{u}\geq \|u\|$ for each $u\in S$,    an appeal to  \cite[Proposition 2.1]{Cross} yields that $S$ is an operator range in $E$.  Moreover,  because of the definition of $M$ we have $M\subset S$,  and thanks to inequality (6) from the previous step we get $R \subset S$.  From the latter it follows that $S$ is dense in $E$. On the other hand, by property (1) from Step 1, 
		$$\|n f_n\| \geq n f_n (w_n) = n \quad \text{for \ each} \quad n \geq 1,$$
		hence the sequence $\{n f_n\}_n$ is not bounded,  and thanks to the uniform boundedness principle it follows that $S$ is not closed.  In particular,  $S \in \mathcal R_d (E)$,  and thus $\cod_E (S) = \infty$. Since $R + M \subset S$, it follows that $\cod (R + M) = \infty$.\par

		It remains to prove that $T$ has the Kato property for $M$.   We claim that the biorthogonal system $\{z_n,f_n\}_n\subset E\times M^{\perp}$ satisfies 
the hypothesis of Corollary \ref{cr:sw} with $\alpha= 3$. 	Since 
		\begin{equation}\label{th:fltz:eq:7}
			\|Q_{M} (z_n)\| = \max \{f (z_n) : f \in S_{M^{\perp}}\} \quad \mathrm{for \ each} \quad n \geq 1,
		\end{equation}
		we get
		$$\|Q_{M} (z_n)\| \geq \|f_n\|^{-1} f_n (z_n) \quad \mathrm{for \ each} \quad n \geq 1.$$
		Therefore, according to property (3) from Step 1 we have 
			\begin{equation}\label{th:fltz:eq:cota1}
			\|T (z_n)\| \leq 3^{-2n} \|Q_{M} (z_n)\| \quad \mathrm{for \ each} \quad n \geq 1.
			\end{equation}
		Now, we will show that
		 	\begin{equation}\label{th:fltz:eq:cota2}
				\|Q_{M} (z_n)\| \leq 3^n \|f_n\|^{-1} f_n (z_n) \quad \mathrm{for \ each} \quad n \geq 1.
			\end{equation}
		First, let us see that
		\begin{equation}\label{th:fltz:eq:9}	
			M^{\perp} = [f_1, ..., f_k] \oplus \overline{[\{f_i : i > k\}]}^{w^*} \quad \mathrm{for \ each} \quad k \geq 1.
		\end{equation}
		Fix $k \geq 1$. Let $h \in M^{\perp}$. Since $\{f_i\}_i$ is linearly weak$^*$-dense in $(E/M)^*$, there exists a net $\{\phi_{\alpha}\}_{\alpha \in \Lambda} \subset \spn \{f_i : i \geq 1\}$ such that $\phi_{\alpha} \xrightarrow{w^*} h$.  For each $\alpha\in \Lambda$ there exists $\{\lambda_{\alpha,i}\}_i \in c_{00}$ such that
		$$\phi_{\alpha} = \sum_i \lambda_{\alpha,i} f_i.$$
		Taking into account that $f_j (w_i) = \delta_{i,j}$ for $i, j \geq 1$ we get $ \phi_{\alpha} (w_i) = \lambda_{\alpha,i}$ for all $\alpha \in \Lambda$ and $i\geq 1$.  Since $\phi_{\alpha} \xrightarrow{w^*} h$ we obtain
		$$\lim_{\alpha} \lambda_{\alpha,i} = \lim_{\alpha} \phi_{\alpha} (w_i) = h (w_i)$$
		for every $i \geq 1$.  For each $\alpha \in \Lambda$ we write 
		$$\widetilde{\phi}_{\alpha} = \sum_{i > k} \lambda_{\alpha,i} f_i.$$
		It is clear that
		$$\widetilde{\phi}_{\alpha} \xrightarrow{w^*} h - (h (w_1) f_1 + \cdots + h (w_k) f_k).$$
		Bearing in mind that $\{\widetilde{\phi}_{\alpha}\}_{\alpha \in \Lambda} \subset \spn \{f_i : i > k\}$, we deduce that
		$$h - (h (w_1) f_1 + \cdots + h (w_k) f_k) \in \overline{[\{f_i : i > k\}]}^{w^*},$$
		which implies \eqref{th:fltz:eq:9}.\par
		Now,  we shall show that 
		\begin{equation}\label{th:fltz:eq:10}
			\overline{[\{f_i : i > k\}]}^{w^*} = [Q_M (v_1), ..., Q_M (v_k)]^{\perp} \quad \mathrm{for \ each} \quad k \geq 1.
		\end{equation}
		Identifying the functionals $f_n$ with elements from $(E/M)^*$, thanks to property (4) from Step 1 we have $\{f_i : i > k\} \subset [Q_M (v_1), ..., Q_M (v_k)]^{\perp}$. Hence
		\begin{equation}\label{auxiliar-step2}
			\overline{[\{f_i : i > k\}]}^{w^*} \subset [Q_M (v_1), ..., Q_M (v_k)]^{\perp}.
		\end{equation}
		On the other hand, it is clear that
		$$\cod_{M^{\perp}} \, \left( [Q_M (v_1), ..., Q_M (v_k)]^{\perp} \right) = k,$$
		and because of \eqref{th:fltz:eq:9}, the subspace $\overline{[\{f_i : i > k\}]}^{w^*}$ also has codimension $k$ in $M^{\perp}$. From this and \eqref{auxiliar-step2} it follows \eqref{th:fltz:eq:10}.\par
		Now,  fix $n \geq 1$,  and pick $f \in S_{M^{\perp}}$ such that $f (z_n) = \|Q_M (z_n)\|$. Because of \eqref{th:fltz:eq:9}, there exist $\lambda_1, ..., \lambda_n \in \R$ and $g \in \overline{[\{f_i : i > n\}]}^{w^*}$ such that
		$$f = \lambda_1 f_1 + \ldots + \lambda_n f_n + g.$$
		Bearing in mind that $\lambda_2 f_2 + \ldots + \lambda_n f_n + g \in \overline{[\{f_i: i > 1\}]}^{w^*}$,  thanks to \eqref{th:fltz:eq:10} we get
		\begin{equation*}
			1 = \|f\| \geq |f (v_1)| = |\lambda_1| |f_1 (v_1)|.
		\end{equation*}
		Since, by property (4) from Step 1 we have $f_1 (v_1) \geq 2^{-1} \|f_1\|$,   we obtain $1 \geq 2^{-1} |\lambda_1| \|f_1\|$,  that is,
		$$|\lambda_1| \leq 2 \|f_1\|^{-1}.$$
		Using subsequently \eqref{th:fltz:eq:10}, property (5) from Step 1 and the former inequality,  it follows that
		\begin{align*}
			1 & =  \|f\| \geq |f (v_2)| = \big{|} \lambda_2 f_2 (v_2) + \lambda_1 f_1 (v_2) \big{|} \geq |\lambda_2| |f_2 (v_2)| - |\lambda_1| |f_1 (v_2)| \\
			& \geq 2^{-1} |\lambda_2| \|f_2\| - 2.
		\end{align*}
		Therefore, 
		$$|\lambda_2| \leq 6 \|f_2\|^{-1}.$$
		Repeating this process $n$ times, we get
		$$|\lambda_n| \leq 2 \cdot 3^{n-1} \|f_n\|^{-1}.$$
		Taking into account property (2) from Step 1, we have $f_1 (z_n) = \cdots = f_{n-1} (z_n) = g (z_n) = 0$. This and the previous inequality imply 
		\begin{equation}\label{th:fltz:eq:11}
			f (z_n) = \lambda_n f_n (z_n) \leq 2 \cdot 3^{n-1} \|f_n\|^{-1} f_n (z_n),\nonumber
		\end{equation}
		which, together with \eqref{th:fltz:eq:7}, yields \eqref{th:fltz:eq:cota2}. Since the biorthogonal system $\{z_n, f_n\}_n$ satisfies \eqref{th:fltz:eq:cota1} and \eqref{th:fltz:eq:cota2},  Corollary \ref{cr:sw} guarantees that $T$ has the Kato property for $M$, and the proof of Step 2 is completed.
	
	\textsc{Step 3}: The statement of the theorem holds true if $L = \{0\}$.
	
	We may assume that $\overline{T(E)} = F$. Since $T$ has the Kato property, by Step 2 there is a closed subspace $M \subset E$ such that $E/M$ is separable, $\cod_E (R + M) = \infty$ and $T$ has the Kato property for $M$. Thanks to Proposition \ref{pr:kato}  it follows that the induced operator $\widetilde{T}_M : E/M \rightarrow F/\overline{T(M)}$ defined as $\widetilde{T}_M (Q_M (u)) = Q_{\overline{T(M)}} (T (u))\, $  $(u \in E)$, has the Kato property.  Set $\widetilde{R} = Q_M (R)$.  Since $\cod_E \, (M + R) = \infty$ we have $M + R \in \mathcal R_d (E)$,  therefore $\widetilde{R} \in \mathcal R_d (E/M)$ (in particular, $\widetilde{R}$ is infinite-codimensional in $E/M$).  Since $E/M$ is separable,  we can apply Theorem \ref{main-FLTZ} to the operator $\widetilde{T}_M$,  the trivial subspace of $E/M$ and the operator range $\widetilde{R}$.  Therefore,  there exists a closed subspace $\widetilde{E}_1 \subset E/M$ such that $\widetilde{E}_1 \cap \widetilde{R} = \{0\}$ and
	$\overline{\widetilde{T}_M (\widetilde{E}_1)} = \overline{\widetilde{T}_M (E/M)}.$
	Let us write 
	$$E_1 = Q_M^{-1} (\widetilde{E}_1).$$
	It is clear that $M \subset E_1$, $\dim \, (E_1/M) = \infty$, $Q_M (E_1) \cap \widetilde{R} = \{0\}$ and
$\overline{\widetilde{T}_M (Q_M (E_1))} = \overline{\widetilde{T}_M (E/M)}.$ Since $Q_M (E_1) \cap \widetilde{R} = \{0\}$,  thanks to  \cite[Theorem 2.4]{Cross} it follows that $\cod_{E/M} \, (Q_M (E_1) + \widetilde{R}) = \infty$,  and thus $\cod_E \, (E_1 + R) = \infty$.  To finish,  we shall show that 
	$$\overline{T (E_1)} = \overline{T (E)}.$$
	Pick $f \in F^*$ such that $f|_{\overline{T(E_1)}} = 0$. Then $f \in \overline{T (E_1)}^{\perp}$, in particular $f \in \overline{T (M)}^{\perp}$,  so $f$ can be identified with a functional on $F/\overline{T(M)}$ which vanishes on $\overline{T(E_1)}/\overline{T(M)}$.  Bearing in mind that 
	\begin{align*}
		\overline{Q_{\overline{T (M)}} (T (E_1))} & = \overline{\widetilde{T}_M (\widetilde{E}_1)} = \overline{\widetilde{T}_M (E/M)} = \overline{Q_{\overline{T (M)}} (T (E))} \\ &= F/\overline{T(M)}
	\end{align*}
	it follows that $f = 0$.
		
	\textsc{Step 4}: The statement holds true for every closed subspace $L \subset E$.\par
		As before, we can assume that $\overline{T(E)} = F$. Using again Proposition \ref{pr:kato}, it follows that the induced operator $\widetilde{T}_L : E/L \rightarrow F/\overline{T(L)}$ has the Kato property.  On the other hand,  as $\cod_E \, (L + R) = \infty$, we deduce that $Q_L (R) \in \mathcal R_d (E/L)$.  Applying the previous step to the operator $\widetilde{T}_L$ and the operator range $Q_L (R)$ we deduce the existence of a closed subspace $\widetilde{E}_1 \subset E/L$ such that $(E/L)/\widetilde{E}_1$ is separable and infinite-dimensional, $\cod_{E/L} (\widetilde{E}_1 + Q_L (R)) = \infty$ and 
	$\overline{\widetilde{T}_L (\widetilde{E}_1)} = \overline{\widetilde{T}_L (E/L)}.$
	Set $E_1 = Q_L^{-1} (\widetilde{E}_1)$. It is immediate that $E/E_1$ is separable and infinite-dimensional and $\cod_E \, (E_1 + R) = \infty$. The fact that
	$\overline{T (E_1)} = \overline{T (E)}$ follows arguing as in the previous step.	
			\hfill $\square$

In the rest of the section we give some applications of Theorem \ref{th:main}.  The first one constitutes  a characterization of dense-range operators that admit  dense-range restrictions to infinite-codimensional subspaces of their domain.  

\begin{cor}
Let $T:E\to F$ be a dense-range operator between Banach spaces.  If $E$ has a separable quotient,  then $T$ has the Kato property if,  and only if,   there exists a closed subspace $E_1\subset E$ such that $\overline{T(E_1)} = F$. 
\end{cor}

\begin{proof}
The ``only if part" follows applying Theorem \ref{th:main} to any proper dense operator range $R$ in $E$ and the subspace $L=\{0\}$. 
Assume there exists a closed (infinite-codimensional) subspace $E_1\subset E$ satisfying $\overline{T(E_1)}=F$.  If $T$ does not have the Kato property then we can find  a closed finite-codimensional subspace $Z\subset E$ such that $T|_Z$ is an isomorphic embedding. Set $X_1= E_1\cap Z$.  Then, $T(X_1)$ is a closed subspace of $F$.  Since $\dim \,  (E_1/X_1) < \infty$ we infer that $T(E_1)$ is closed as well,  therefore $T(E_1)=F$.  But taking into account that $\dim \, (E_1) = \infty$ we get $\dim \, (\ker T) = \infty$,  and consequently $Z\cap \ker T \neq \{0\}$,  which contradicts the injectivity of $T|_Z$. 
\end{proof}

As we mentioned,  Johnson \cite{Johnson} (see also \cite{James}) proved that if $X$ and $Y$ are proper quasicomplementary subspaces of a Banach space $E$,  and $X$ is weakly compactly generated, then there is a closed subspace $X_1\subset X$ such that $\dim \, (X/X_1)= \infty$ and $X_1$ is a quasicomplement of $Y$.  The following corollary of Theorem \ref{th:main} provides a generalization of this result under the only assumption that the subspace $X$ has a separable quotient. 

\begin{cor}\label{cr:james}
	Let $X$ and $Y$ be closed subspaces of a Banach space $E$ such that $X+Y$ is dense in $E$.  If $X$ has a separable quotient and $X+Y$ is not closed  then there exists a closed subspace $X_1 \subset X$ such that $\dim \, (X/X_1)= \infty$ and $X_1 + Y$ is dense in $E$.  
\end{cor}

\begin{proof}
	Consider the operator $T = Q_Y|_X:X\to E/Y$.  Since $X+Y$ is dense in $E$ we have that $\overline{T(X)}=E/Y$.  Moreover, as $X+Y$ is not closed,   thanks to Corollary \ref{cr:katoq1} it follows that $T$ has the Kato property.  Using Theorem \ref{th:main} we deduce the existence of a closed subspace $X_1 \subset X$ such that $X/X_1$ is infinite-dimensional and separable, and $\overline{T (X_1)}  = E/Y$.
	The latter yields that $X_1 + Y$ is dense in $E$.
\end{proof}

Next,  we provide a result on spaceability in dual spaces.  In the paper \cite{P81}, Plichko proved that if $E$ is a separable Banach space then, for every Banach space $F$ and every operator $T : E \rightarrow F$ such that $T(E)$ is proper and dense in $F$, there exists a closed subspace $X \subset E$ such that $T^* (F^*) \cap X^{\perp} = \{0\}$.  The next result gives an extension of this statement assuming only that $E$ has a separable quotient.

\begin{cor}\label{cr:plichko}
	Let $T : E \rightarrow F$ be a dense-range operator between Banach spaces, where $E$ has a separable quotient, and let $R \in \mathcal R_d (E)$. If $T(E)$ is not closed then there exists a closed subspace $X \subset E$ such that:
	\begin{enumerate}
	\item[(a)] $E/X$ is separable and infinite-dimensional, 
	\item[(b)]	$\cod_E \, (R + X) = \infty$, and
	\item[(c)]	$T^* (F^*) \cap X^{\perp} = \{0\}.$
	\end{enumerate}
\end{cor}

\begin{proof}
	Since $T(E)$ is dense and not closed in $F$, it follows that $T$ has the Kato property. Therefore, according to Theorem \ref{th:main} there exists a closed subspace $X \subset E$ such that $E/X$ is separable and infinite-dimensional,  $\cod_E \, (R + X) = \infty$ and $\overline{T (X)} = F$. Pick $f \in F^*$ such that $T^* (f) \in X^{\perp}$. Then,  for every $x\in X$ we have $$\langle f , \, T(x) \rangle = \langle T^* (f), \,  x \rangle = 0 .$$ Since $\overline{T (X)}=F$ it follows that $f = 0$.  Consequently, $T^* (F^*) \cap X^{\perp} = \{0\}$.
\end{proof}

We end this section with a variant of the previous result. 
\begin{cor}
Let $T : E \to F$ be an operator with non-closed range between Banach spaces. If $E$ has a separable quotient then, for each closed subspace $M \subset T (F)$ there exists a closed subspace $X\subset E$ such that:  
\begin{enumerate}
\item $E/X$ is separable and infinite-dimensional,
\item $M \subset T (X)$,  and  
\item $T^*(F^*) \cap X^{\perp} = \{0\}$. 
\end{enumerate}
\end{cor}

\begin{proof}
	Set $L = T^{-1}(M)$. By Corollary \ref{cr:kato1}, $T$ has the Kato property for $L$. Therefore, according to Theorem \ref{th:main} there exists a closed subspace $X \subset E$ such that 
$E/X$ is infinite-dimensional and separable,  $L \subset X$ (hence $M\subset T(X)$),   and $\overline{T(X)} =F$ (consequently  $T^*(F^*)\cap X^{\perp} = \{0\}$). 
\end{proof}

\section{Operators with the Kato property in Banach spaces with weak$^*$-separable dual}

In this section,  we provide the next extension of Theorem \ref{main-FLTZ}. 

\begin{thm} \label{th:main-ws}
	Let $E$ and $F$ be Banach spaces,  let $L \subset E$ be a closed subspace,  let $T : E \rightarrow F$ be an operator with the Kato property for $L$,  and $R \in \mathcal R_d (E)$ such that $$R\cap L = \{0\}.$$ If $(E/L)^*$ is weak$^*$-separable then there exists a closed subspace $E_1 \subset E$ such that 
	\begin{enumerate}
\item[(a)] $L \subset E_1$, 
\item[(b)] $E/E_1$ is separable and infinite-dimensional, 
\item[(c)] $E_1\cap R = \{0\}$, and
\item[(d)] $\overline{T (E_1)} = \overline{T (E)}.$
\end{enumerate}
	
\end{thm}

In the proof of this theorem, we will make use of the following lemma.

\begin{lema}[{\cite[Lemma 3]{Fonf}, see also \cite[Lemma 2.2]{JL1}}]\label{lm:jl1}
	Let $E$ be a Banach space and let $R \in \mathcal R_d (E)$. Then, for every pair of sequences $\{v_n\}_n \subset S_E$ and $\{\varepsilon_n\}_n \subset (0, \infty)$ there exist sequences $\{w_n\}_n \subset S_E$ and $\{\gamma_n\}_n \subset (0, 1)$ such that $\sum_n \gamma_n \leq 1$,
	$$\spn \left( \clco \{\pm \gamma_n w_n\}_n\right) \cap R = \{0\},$$
	and
	$$\|v_n - w_n\| \leq \varepsilon_n \quad \text{for \ each} \quad n \geq 1.$$
\end{lema}

\begin{proof}[Proof of Theorem \ref{th:main-ws}]
Assume that $\overline{T (E)} = F$. First,  we prove the result in the case that $L = \{0\}$. According to Step 2 from the proof of Theorem \ref{th:main} there exists a closed subspace $M \subset E$ such that $E/M$ is separable, $\cod_E \, (M + R) = \infty$ and $T$ has the Kato property for $M$. Thus,  according to Proposition \ref{pr:kato} there exists a sequence $\{u_n\}_n \subset E$ such that $\{Q_M (u_n)\}_n$ is a normalized basic sequence in $E/M$ and $\|T (u_n)\| \leq 2^{-n-2}$ for each $n \geq 1$. Denote by $\beta$ the basic constant of $\{Q_M (u_n)\}_n$.\par
	Notice that $M + R \in \mathcal R_d (E)$. Applying Lemma \ref{lm:jl1} we can find sequences $\{u_n'\}_n \subset S_E$ and $\{\gamma_n\}_n \subset (0, \infty)$ such that $\sum_n \gamma_n \leq 1$,
	\begin{equation}\label{th:main-ws:eq:4}
		\spn \left( \clco \{\pm \gamma_n u_n'\}_n \right) \cap (M + R) = \{0\},
	\end{equation}
	and
	\begin{equation*}
		\left \|u_n' - \frac{u_n}{\|u_n\|} \right \| \leq 2^{-n-2} (2 + \|T\|)^{-1} \|u_n\|^{-2} \beta^{-1} \quad \mathrm{for \ each} \quad n \geq 1.
	\end{equation*}
	
	Dividing by 2 if needed, we may assume that $\sum_n \gamma_n \leq 2^{-1}$.  For each $n\geq 1$ we set $u_n'' = \|u_n\| u_n'$. Then
	\begin{equation}\label{th:main-ws:eq:1}
		\|u_n'' - u_n\| \leq 2^{-n-2} (2 + \|T\|)^{-1} \|u_n\|^{-1} \beta^{-1} \quad \mathrm{for \ each} \quad n \geq 1,
	\end{equation}
	Let us write 
	$$w_n = \frac{u_n''}{\|Q_M (u_n'')\|},  \quad n \geq 1.$$
	Fix $n \geq 1$.  Because of \eqref{th:main-ws:eq:1} we have
	\begin{align*}
		\|Q_M (u_n'')\| & \geq \|Q_M (u_n)\| - \| Q_M (u_n'' - u_n)\| \geq 1 - \|u_n'' - u_n\|  \\
		& \geq 1 - 2^{-n-2} (2 + \|T\|)^{-1} \|u_n\|^{-1} \beta^{-1}, 
	\end{align*}
	and 
\begin{equation}	
	\|Q_M (u_n'')\| \leq \|Q_M (u_n)\| + \| Q_M (u_n'' - u_n)\| \leq 1 + 2^{-n-2} (2 + \|T\|)^{-1} \|u_n\|^{-1} \beta^{-1}.
	\end{equation}
	From the first of the two former inequalities we get 
	\begin{align*}
		\|Q_M (u_n'')\|^{-1} \leq & \ [1 - 2^{-n-2} (2 + \|T\|)^{-1} \|u_n\|^{-1} \beta^{-1}]^{-1}  \\
		\leq & \ 1 + 2^{-n-2} (2 + \|T\|)^{-1} \|u_n\|^{-1} \beta^{-1} [1 - 2^{-n-2} (2 + \|T\|)^{-1} \|u_n\|^{-1} \beta^{-1}]^{-1}  \\
		\leq & \ 1 + 2^{-n-1} (2 + \|T\|)^{-1} \|u_n\|^{-1} \beta^{-1},
	\end{align*}
	and the second one yields 
	$$\|Q_M (u_n'')\|^{-1} \geq 1 - 2^{-n-1} (2 + \|T\|)^{-1} \|u_n\|^{-1} \beta^{-1}.$$
	Since $\|u_n\| \geq 1$, we get
	\begin{align*}
		\|Q_M (w_n - u_n)\| & = \left \| \|Q_M (u_n'')\|^{-1} Q_M (u_n'' - u_n) + (\|Q_M (u_n'')\|^{-1} - 1) Q_M (u_n)  \right \|  \\
		& \leq \|Q_M (u_n'')\|^{-1} \|Q_M (u_n'' - u_n)\| + \left | \|Q_M (u_n'')\|^{-1} - 1 \right | \|Q_M (u_n)\|  \\
		& \leq \|Q_M (u_n'')\|^{-1} \|u_n'' - u_n\| + \left | \|Q_M (u_n'')\|^{-1} - 1 \right | \|u_n\|  \\
		& \leq 2^{-n-1} (2 + \|T\|)^{-1} \beta^{-1} + 2^{-n-1} (2 + \|T\|)^{-1} \beta^{-1}  \\
		& \leq 2^{-n-1} \beta^{-1}.
	\end{align*}
	According to the stability theorem (c.f. \cite[Theorem 4.23]{Fabian}) it follows that $\{Q_M (w_n)\}_n$ is a normalized basic sequence in $E/M$. Moreover, for each $n \geq 1$ we have 
	\begin{align*}
		\|T (w_n)\| & \leq \|Q_M (u_n'')\|^{-1} \|T (u_n'')\| \leq 2 \|T (u_n'')\| \leq 2 \|T (u_n)\| + 2 \|T (u_n'' - u_n)\|  \\
		& \leq 2 \cdot 2^{-n-2} + 2 \cdot 2^{-n-2} \leq 2^{-n}.
	\end{align*}

Because of the weak$^*$-separability of $E^*$,  there exists a sequence $\{e_n^*\}_n\subset S_{E^*}$ which is total over $E$. Let $n_0 \geq 1$ such that the basic constant of $\{w_n\}_n$ is smaller than $2^{n_0/2}$, let $\{w_n^*\}_n \subset M^{\perp}$ be the sequence of biorthogonal functionals of $\{Q_M (w_n)\}_n$, and set
	$$Y_n = \bigcap_{i=1}^{n} \ker e_i^* \quad \text{for each} \quad n \geq 1.$$
	Applying Lemma \ref{lm:sw} to the biorthogonal system $\{w_n, w_n^*\}_{n \geq n_0}$ and the sequence of subspaces $\{Y_n\}_n$ (with $\alpha=\sqrt{2}$) we deduce the existence of a block sequence $\{v_n\}_n$ of $\{w_n\}_n$ such that, for every $n \geq 1$:
	\begin{itemize}
		\item[(1)] $v_n \in Y_n$, 
		\item[(2)] $\|Q_M (v_n)\| = 1$, and
		\item[(3)] $\|T (v_n)\| \leq 2^{-n-2}$.
	\end{itemize}
	
	Since $\{v_n\}_n$ is a block sequence of $\{w_n\}_n$ and for each $n \geq 1$, $w_n$ and $u_n'$ are collinear, we have that $\{Q_M (v_n)\}_n$ is a block sequence of $\{Q_M (u_n')\}_n$. Taking into account that  $\{Q_M (u_n')\}_n$ is basic, it follows that $\{Q_M (v_n)\}_n$ is basic as well. In addition, since $\{v_n\}_n$ is a block sequence of $\{u_n'\}_n$, for each $n \geq 1$ there exist $p_n, q_n \geq 1$ and $\nu_{p_n+1}, ..., \nu_{p_n+q_n} \in \R$ such that
	$$v_n = \nu_{p_n+1} u_{p_n+1}' + \cdots + \nu_{p_n+q_n} u_{p_n+q_n}'.$$
	Let us write  $\sigma_n = \sum_{i=1}^{q_n} |\nu_{p_n+i}|$, $\xi_n = \min \{\gamma_{p_n+i} : \, 1 \leq i \leq q_n\}$, and $\mu_n = q_n^{-1} \sigma_n^{-1} \xi_n$ for each $n \geq 1$. We claim that
	\begin{equation}\label{th:main-ws:eq:2}
		\spn \left( \clco (\{\pm \mu_n v_n\}_n) \right) \cap (M + R) = \{0\}.
	\end{equation}
	Indeed, bearing in mind that $\{u_n'\}_n \subset S_E$, $\sigma_n \geq \|v_n\|$ for all $n \geq 1$, and $\{\gamma_n\}_n \in \ell_1$  it follows that $\sum_n \mu_n \|v_n\| < \infty$, hence $\clco (\{\pm \mu_n v_n\}_n)$ coincides with the set of sums $\sum_n \lambda_n \mu_n v_n$, where $\{\lambda_n\}_n \subset [-1,1]$ is a sequence such that $\sum_n |\lambda_n| \leq 1$. Let $\{\lambda_n\}_n \subset [-1,1]$ be such a sequence. Then
	$$\sum_n \lambda_n \mu_n v_n = \sum_n q_n^{-1} \lambda_n \sigma_n^{-1} \xi_n (\nu_{p_n+1} u_{p_n+1}' + \cdots + \nu_{p_n+q_n} u_{p_n+q_n}').$$
	Since $\sigma_n^{-1} \xi_n \nu_{p_n+i} \leq \gamma_{p_n+i}$ for all $1 \leq i \leq q_n$ and $n\geq 1$,   we have that $\sum_n \lambda_n \mu_n v_n \in \clco (\{\pm \gamma_k u_k'\}_k)$. Thus, by \eqref{th:main-ws:eq:4}, we obtain \eqref{th:main-ws:eq:2}.\par		
	Set $\alpha_1 = 3^{-1}\mu_1 \|v_1\|^{-1}$,  and for $n \geq 2$,
	$$\alpha_n = \min \left \{3^{-n} \mu_n \|v_n\|^{-1}, \,\,  \min_{1 \leq i \leq n-1} \{ |e_n^* (v_i)|^{-1} : \, \,  e_n^* (v_i) \neq 0\} \right\},$$
	and define,  for each $u\in E$, 
	$$\varphi (u) = u + \sum_n \alpha_n e_n^* (u) v_n. $$
	It is clear that $\varphi$ is an isomorphism with 
	$$\|\varphi (u)\| \geq 2^{-1} \|u\|, \quad u \in E,$$
	hence $\|\varphi^{-1}\| \leq 2$.  Next,  we prove the following claim.

	\textsc{Claim:}  The operator $T$ satisfies the Kato property for the subspace $\varphi(M)$ and 
	\begin{equation}	 
	 \varphi (M + R) \cap (M + R) = \{0\}.\label{main2-aux}
	 \end{equation}
	First, for $u \in E$ we will calculate a lower bound of $\|Q_{\varphi(M)} (\varphi (u))\|$. If $u \in E$ then
	\begin{align*}
		\|Q_{\varphi(M)} (\varphi (u))\| & = \inf \{ \|\varphi (u) - \varphi (v)\| : v \in M \} \geq \|\varphi^{-1}\|^{-1} \inf \{\|u - v\| : v \in M\} \\
		& = \|\varphi^{-1}\|^{-1} \|Q_M (u)\| \geq 2^{-1} \|Q_M (u)\|.
	\end{align*}
	Now, for each $n \geq 1$, we will find an upper bound for $\|T (\varphi (v_n))\|$. Fix $n \geq 1$. Since $v_n \in \bigcap_{i=1}^n \ker e_i^*$ we have
	$$\varphi (v_n) = v_n + \sum_k \alpha_k e_k^* (v_n) v_k = v_n + \sum_k \alpha_{n+k} e_{n+k}^* (v_n) v_{n+k}.$$
	On the other hand, as $\|T (v_k)\| \leq 2^{-k-2}$ for each $k \geq 1$, we have 
	\begin{align*}
		\|T (\varphi (v_n))\| &\leq \|T (v_n)\| + \sum_k \alpha_{n+k} |e_{n+k}^* (v_n)| \|T (v_{n+k})\|\\ &\leq 2^{-n-2} \left(1 + \sum_k 2^{-k} \alpha_{n+k} |e_{n+k}^* (v_n)| \right).
	\end{align*}
	Since $\alpha_{n+k} \leq |e_{n+k}^* (v_n)|^{-1}$ for each $k \geq 1$ such that $e_{n+k}^* (v_n) \neq 0$, we get
	$$\sum_k 2^{-k} \alpha_{n+k} |e_{n+k}^* (v_n)| \leq 1,$$
	and consequently 
	\begin{equation}\label{th:main-ws:eq:3}
	 	\|T (\varphi (v_n))\| \leq 2^{-n-1} = 2^{-n-1} \|Q_M (v_n)\| \leq 2^{-n}\|Q_{\varphi(M)} (\varphi (v_n))\|.
	\end{equation}
	
	Bearing in mind that $\varphi$ is an isomorphism we have that the operator $\Phi : E/M \rightarrow E/\varphi(M)$ defined by the formula $\Phi (Q_M (u)) = Q_{\varphi(M)} (\varphi (u))$ for $u \in E$ is also an isomophism. Therefore, since $\{Q_M (v_n)\}_n$ is a basic sequence in $E/M$ it follows that $\{Q_{\varphi(M)} (\varphi (v_n))\}_n$ is a basic sequence in $E/\varphi(M)$. Thus, and taking into account \eqref{th:main-ws:eq:3},  thanks to Proposition \ref{pr:kato} it follows $T$ satisfies the Kato property for $\varphi (M)$.  Now,  let us check \eqref{main2-aux}.  Pick $u \in M + R$ such that $\varphi (u) \in M + R$. Then
	 $$\sum_n \alpha_n e_n^* (u) v_n = \varphi (u) - u \in M + R.$$
	 Since $\{\alpha_n \mu_n^{-1} e_n^* (u)\}_n \in \ell_1$, we get
	 $$\sum_n \alpha_n e_n^* (u) v_n \in \spn \left(\clco \{\pm \mu_n v_n\}_n \right).$$
	 Therefore, by \eqref{th:main-ws:eq:2}, it follows that $\sum_n \alpha_n e_n^* (u) v_n = 0$. Since $\{v_n\}_n$ is minimal, we have that $e_n^* (u) = 0$ for each $n \geq 1$. Taking into account that $\{e_n^*\}_n$ is total over $E$, it follows that $u = 0$.  Thus,  $\varphi(M+R) \cap (M+R) = \{0\}$,  and the claim is proved.  \par
	 
	 Next,  we will build the subspace $E_1$ with the required properties.  According to Proposition \ref{pr:kato},  the operator $\widetilde{T}_{\varphi(M)} : E/\varphi(M) \rightarrow F/\overline{T(\varphi(M))}$ defined by the formula $\widetilde{T}_{\varphi(M)} (Q_{\varphi(M)} (u)) = Q_{\overline{T({\varphi(M)})}} (T (u))\, $  ($u \in E$) has the Kato property. Set
	 $\widetilde{R} = Q_{\varphi(M)} (R).$
	 Since ${\varphi(M)} \cap R = \{0\}$  we get $\cod ({\varphi(M)} + R) = \infty$, hence ${\varphi(M)} + R \in \mathcal R_d (E)$. Thus $\widetilde{R} \in \mathcal R_d (E/{\varphi(M)})$.  Because of the separability of $E/{\varphi(M)}$,  we can apply Theorem \ref{main-FLTZ} to obtain  a closed infinite-dimensional subspace $\widetilde{E}_1 \subset E/{\varphi(M)}$ such that $\widetilde{E}_1 \cap \widetilde{R} = \{0\}$ and
	$\overline{\widetilde{T}_{\varphi(M)} (\widetilde{E}_1)} = \overline{\widetilde{T}_{\varphi(M)} (E/{\varphi(M)})}.$
	Let us write $$E_1 = Q_{\varphi(M)}^{-1} (\widetilde{E}_1).$$ Then ${\varphi(M)} \subset E_1$, $Q_{\varphi(M)} (E_1) \cap \widetilde{R} = \{0\}$ and
		$$\overline{\widetilde{T}_{\varphi(M)} (Q_{\varphi(M)} (E_1))} = \overline{\widetilde{T}_{\varphi(M)} (E/{\varphi(M)})}.$$
	Let us show that $\overline{T (E_1)} = \overline{T (E)}$. Pick $f \in F^*$ such that $f|_{\overline{T (E_1)}} = 0$. Then $f$ identifies with a functional on $F/\overline{T(M)}$ which vanishes on $\overline{T(E_1)}/\overline{T(M)}$. Bearing in mind that
	\begin{align*}
		\overline{Q_{\overline{T(M)}} (T (E_1))} & = \overline{\widetilde{T}_{\varphi(M)} (Q_{\varphi(M)} (E_1))} 
		= \overline{\widetilde{T}_{\varphi(M)} (E/{\varphi(M)})} \\ 
		& = \overline{Q_{\overline{T(M)}} (T (E))} 
		= F/\overline{T(M)},
	\end{align*}
	we conclude that $f = 0$, and the theorem is proved whenever $L=\{0\}$.
	
To finish,  suppose that $L$ is a non-trivial closed subspace of $E$ and that $(E/L)^*$ is weak$^*$-separable.   As the operator $T$ has the Kato property for $L$,  we have that the induced operator $\widetilde{T}_L : E/L \rightarrow F/\overline{T(L)}$ has the Kato property.  According to the particular case,  applied to the space $E/L$,  the operator  $\widetilde{T}_L$ and the (proper dense) operator range $Q_L(R)$,  we deduce the existence of an infinite-dimensional closed subspace $\widetilde{E}_1 \subset E/L$ such that
	$$\overline{\widetilde{T}_L (\widetilde{E}_1)} = \overline{\widetilde{T}_L (E/L)} \quad \text{and} \quad \widetilde{E}_1 \cap Q_L(R) = \{0\}.$$
	Now, set $E_1 = Q_L^{-1} (\widetilde{E}_1)$. Then $L \subset E_1$ and $\dim \, (E_1/L) = \infty$, $E/E_1$ is separable and infinite-dimensional,  
	$$\overline{\widetilde{T}_L (E_1/L)} = \overline{\widetilde{T}_L (E/L)} \quad \text{and} \quad (E_1/L) \cap Q_L(R) = \{0\}.$$
	From the first of the last two properties, arguing as in the particular case we get $\overline{T(E_1)} =\overline{T(E)}.$  On the other hand,  if $x\in R\cap E_1$ then $Q_L(x)\in (E_1/L) \cap Q_L(R)$,  hence (according to the last property), $x\in L$,  and taking into account that $R\cap L= \{0\}$ it follows that $x=0.$ Consequently,  $E_1\cap R = \{0\}$. 
\end{proof}

\begin{rem}\normalfont
Drewnowski \cite{D} proved that any infinite-dimensional Banach space contains a proper dense linear subspace $V$ such that,  if $X$ is a closed subspace of $E$ with $X\cap V = \{0\}$,  then $X$ is finite-dimensional.  Therefore,  it is essential to require that the subspace $R$ in the previous theorem is a proper dense operator range (or a non-barrelled proper dense subspace) of $E$.
\end{rem}

Now,  we give some applications of the former theorem.  As an easy consequence we obtain the following result. 

\begin{cor}\label{cr:1}
	Let $L$ be a closed subspace of a Banach space $E$ and let $R \in \mathcal R_d (E)$ such that $R \cap L = \{0\}$. If $L^{\perp}$ is weak$^*$-separable then, there exists a closed subspace $E_1 \subset E$ such that $L \subset E_1$,  $\dim (E_1/L) = \infty$,   $E/E_1$ is separable  and $R \cap E_1 = \{0\}$.
\end{cor}

\begin{proof}
	Set $\widetilde{R} = Q_L (R)$.  Then $\widetilde{R}\in \mathcal R_d (E/L)$.  Since $(E/L)^*$ is weak$^*$-separable,  Theorem \ref{th:main-ws},  applied with the null endomorphism on $E$,  guarantees the existence of a closed subspace $E_1\subset E$ satisfying the specified properties. 
\end{proof}

Next,  we provide some results on quasicomplements.  The first one constitutes  an extension of \cite[Corollary 2.3]{FLTZ}. 

\begin{cor}\label{fltz-1bis}
	Let $X$ be a closed infinite-dimensional subspace of a Banach space $E$ such that $E/X$ has a separable quotient,  let $R\in \mathcal R_d(E)$ with $X\subset R$,  and let $L\subset E$ be a closed subspace such that
\begin{enumerate}
\item[(a)]	
	             $R \cap L = \{0\}$ and 
	             \item[(b)] $(E/L)^*$ is weak$^*$-separable. 
	             \end{enumerate}
	             Then,  there exists a closed subspace $E_1\subset E$ such that 
\begin{enumerate}
	\item $L\subset E_1$, $\dim \, (E_1/L)=\infty$,  $E/E_1$ is separable, 
	\item $R \cap E_1 = \{0\}$,   and 
	\item $E_1+X$ is dense in $E$. 
\end{enumerate}
\end{cor}
\begin{proof}
Let $T=Q_X:E\to E/X$.  Since $\dim (X) = \infty$,   $T$ has the Kato property for any closed subspace of $E$.  An appeal to Theorem \ref{th:main-ws} yields the existence of a closed subspace $E_1\subset E$ satisfying $(1)$,  $(2)$ and $\overline{T(E_1)}=E/X$.  The latter implies $(3)$. 
\end{proof}

Arguing as in the proof of Corollary \ref{cr:james} we obtain the following extension of \cite[Corollary 2.4]{FLTZ}.

\begin{cor}
	Let $E$ be a Banach space. If $X$ and $Y$ are two closed subspaces of $E$  such that $X + Y$ is a proper dense subspace of $E$, $X$ has a separable quotient and $X^*$ is weak$^*$-separable,  then for every $R \in \mathcal R_d (X)$ there exists a closed subspace $X_1 \subset X$ such that $\dim \, (X/X_1)=\infty$, $X_1 + Y$ is dense in $E$, and $R \cap X_1 = \{0\}$.
\end{cor}

A classical theorem of Lindenstrauss and Rosenthal (c.f. \cite[Theorem 5.79]{Hajek},  see also \cite[Corollary 2.9 and Theorem 3.1]{JL2}) asserts that a closed subspace $X$ of a Banach space $E$ is quasicomplemented if $X^*$ is weak$^*$-separable and $E/X$ has a separable quotient.  In view of this result we can wonder if,  for every operator range $R\in \mathcal R_d(E)$ containing $X$,  there exists or not a closed subspace $Y\subset E$ which is a quasicomplement of $X$ satisfying 
$R\cap Y = \{0\}$.  The next result guarantees that the answer to this problem is affirmative if,  and only if,  $E^*$ is weak$^*$-separable. 

\begin{cor}\label{main2}
If $E$ is a Banach space with a separable quotient,  then the following assertions are equivalent: 
\begin{enumerate}
\item $E^*$ is weak$^*$-separable. 
\item For every closed subspace $X\subset E$ such that $E/X$ has a separable quotient and every $R \in \mathcal R_d (E)$ with $X\subset R$ there exists a quasicomplement $Y$ of $X$ such that $R\cap Y = \{0\}$. 
\item For some closed separable subspace $X\subset E$ such that $E/X$ has a separable quotient and every $R \in \mathcal R_d (E)$ with $X\subset R$ there is a quasicomplement $Y$ of $X$ such that and $R\cap Y = \{0\}$.
\end{enumerate}
\end{cor}

\begin{proof} The implication $(1)\Rightarrow (2)$ follows from Corollary \ref{fltz-1bis},  and  $(2)\Rightarrow (3)$ is trivial.  Now,  consider a closed separable subspace $X\subset E$ satisfying $(3)$.  According to \cite[Corollary 2.9]{JL2},  $X$ admits a proper  quasicomplement $Z$ such that 
$E/Z$ is separable.  Set $R = X+Z$.   Then $R\in \mathcal R_d(E)$.  Therefore,  because of the hypothesis,  we can find a new quasicomplement of $X$,  say $Y$,  satisfying
$Y\cap R = \{0\}$.  Since $Y\cap Z = \{0\}$,  we have that the restriction to $Y$ of the quotient map $Q_Z:E\to E/Z$ is one-to-one,  and taking into account that $(E/Z)^*$ is weak$^*$-separable (being the dual of a separable space) it follows that $Y^*$ is weak$^*$-separable as well.  On the other hand,  as $X+Y$ is dense in $E$ we have that $Q_Y(X)$ is dense in $E/Y$,  and because of the separability of $X$ we deduce that $E/Y$ is separable.  Therefore,  both spaces $Y^*$ and $(E/Y)^*$ are weak$^*$-separable.  Taking into account that the property of having weak$^*$-separable dual is a three-space property it follows that  $E^*$ is weak$^*$-separable,  hence $(3)\Rightarrow (1)$.  
\end{proof}

A Banach space $E$ is said to be a \textbf{DENS space} if the density character of $E$ agrees with with the weak$^*$-density character of $E^*$.  The class of such spaces is quite large (it includes,  for instance,  weakly Lindel\"of determined spaces,  c.f.  \cite[Proposition 5.40]{Hajek}).  As an immediate consequence of the previous corollary we obtain the following result. 

\begin{cor}
Let $E$ be a DENS space with a separable quotient.  Then,  $E$ is separable if,  and only if,  for every closed subspace $X\subset E$ and every $R\in \mathcal R_d(E)$ containing $X$ there exists a closed subspace $Y\subset E$ such that $X+Y$ is dense in $E$ and $Y\cap R = \{0\}$. 
\end{cor}

In \cite[Theorem 4.5]{JL3} (see also \cite[Theorem 6.2]{CS}) it is shown that if $E$ is a Banach space with weak$^*$-separable dual,  then for every $R\in \mathcal R_d(E)$ there exist two proper quasicomplementary subspaces $X,Y\subset E$ such that the spaces $E/X$ and $E/Y$ both are separable and  $R\cap (X+Y)= \{0\}$.  The next result provides a strengthening of this fact and \cite[Corollary 2.5]{FLTZ}. 

\begin{cor}\label{cor-2}
Let $E$ be a Banach space,  let $L$ be a closed subspace of $E$ and $R\in \mathcal R_d(E)$ such that $R\cap L = \{0\}$.  If $L^{\perp}$ and $L^*$ are weak$^*$-separable then 
there exist quasicomplementary subspaces $X,Y\subset E$ such that: 
\begin{enumerate}
\item $E/X$ and $E/Y$ are separable,  
\item $L\subset X$,  $\dim (X/L)= \infty$,  and 
\item $R\cap (X+Y)=\{0\}$.
\end{enumerate}
\end{cor}

\begin{proof}
	According to Corollary \ref{cr:1}, there exists a closed subspace $X \subset E$ such that $L \subset X$ and $\dim (X/L) = \infty$ and $R \cap X = \{0\}$. Set $R_1 = R + X$. As $L^*$ and $L^{\perp}$ are weak$^*$-separable, we have that $E^*$ is weak$^*$-separable as well. Since $X \supset R_1$, Corollary \ref{main2} yields the existence of a closed subspace $Y \subset E$ such that $E/Y$ is separable,  $X$ and $Y$ are quasicomplementary and $Y \cap R_1 = \{0\}$.  Pick a vector $u = x + y \in R \cap (X + Y)$. Then $y = u - x \in R_1$,  hence $y = 0$, and so $u = x \in R \cap X$. Since $R \cap X = \{0\}$ we obtain $u = 0$.  Consequently, $R\cap (X + Y)=\{0\}.$
\end{proof}

Plichko \cite{P} (c.f.  \cite[Theorem 5.90]{Hajek}, see also \cite{James,  Johnson}) proved that if $X$ and $Y$ are proper quasicomplementary subspaces of a Banach space $E$,  then $Y$  admits a quasicomplement $Z\subset E$ that contains $X$ as an infinite-codimensional subspace.  The next result yields a refinement of this statement in the case that $E^*$ is weak$^*$-separable.

\begin{cor}\label{cor-3}
	Let $X$ and $Y$ be proper quasicomplementary subspaces of a Banach space $E$.  If $X^*$ and $X^{\perp}$ are weak$^*$-separable then there exists a quasicomplement $Z$ of $Y$ 
such that $Z\supset X$ and $E/Z$ is separable. 
\end{cor}
\begin{proof}
According to \cite[Proposition 4.3 (2)]{JL3},  there exists $R\in \mathcal R_d(E)$ such that $R\supset Y$ and $R\cap X = \{0\}$.  Thus,  thanks to the previous corollary, 
we can find closed quasicomplementary subspaces $Z,W\subset E$ such that $E/Z$ is separable,  $X\subset Z$,  and $R\cap (W+Z) = \{0\}$.  Then,  $Y+Z$ is dense in $E$,  and bearing in mind that $Y\subset R$ we get $Y\cap Z = \{0\}$.  
\end{proof}

We end this paper with a result on spaceability.  As we mentioned, Plichko \cite{P81}  proved that if $E$ is a Banach space with a separable quotient,  then for every operator $T:E\to F$ with dense non-closed range there is a closed infinite-dimensional subspace $X\subset F$ such that $T(E)\cap X = \{0\}$.  Moreover,  Corollary \ref{cr:plichko} ensures the existence of a closed subspace $Y\subset E$ such such that $T^*(F)\cap Y^{\perp} = \{0\}$.  Thus,   it is natural to ask if,  in the case that $F=E$,  
it is or not possible to find a closed infinite-dimensional subspace $X\subset E$ satisfying $$T(E)\cap X = \{0\} \quad \text{and} \quad T^*(E^*)\cap X^{\perp} = \{0\}.$$ The next result guarantees that the answer is affirmative whenever $E^*$ is weak$^*$-separable. 
\begin{cor}\label{plichko-dual-2}
	Let $E$ be a Banach space with a separable quotient and weak$^*$-separable dual, $R \in \mathcal R_d (E)$, $F$ be a Banach space, and $T : E \rightarrow F$ be an operator with proper dense range. Then there exists a closed subspace $X \subset E$ such that:
	\begin{itemize}
		\item[(1)] $E/X$ is separable,
		\item[(2)] $R\cap X=\{0\}$, and
		\item[(3)] $T^* (F^*) \cap X^{\perp}= \{0\}$.
	\end{itemize}
\end{cor}

\begin{proof}
We notice that $T$ has the Kato property.  Thanks to Theorem \ref{th:main-ws} there exists a closed subspace $X\subset E$ such that $E/X$ is separable, $R\cap X = \{0\}$, and $\overline{T(X)}=F$.  From the latter we deduce that $T^* (F^*) \cap X^{\perp} = \{0\}$.
\end{proof}

\noindent {\bf Acknowledgements.}
 M. Jiménez-Sevilla acknowledges the support provided by the grant PID2022-138758NB-I00 from the Ministerio de Ciencia e Innovación of Spain. 
 Sebastián Lajara acknowledges the support provided by the grant  PID2021-122126NB-C32 funded by MCIN/AEI/10.13039/501100011033 and the grant
“ERDF A way of making Europe”, by the European Union.

\end{document}